\documentclass[11pt]{article}
\usepackage[T1]{fontenc}
\usepackage[utf8]{inputenc}
\usepackage{lmodern}
\usepackage{textcomp}
\usepackage[margin=2.3cm]{geometry}
\usepackage{amsmath,amssymb,amsthm}
\usepackage{array,booktabs,longtable}
\usepackage{enumitem}
\usepackage{microtype}
\usepackage{xcolor}
\usepackage[hyphens]{url}
\usepackage[unicode,pdfencoding=auto,colorlinks=true,linkcolor=blue!50!black,urlcolor=blue!50!black,citecolor=green!40!black]{hyperref}
\usepackage{xurl}
\hypersetup{pdftitle={What the Machine Invented: Mining New Mathematical Instruments from OpenAI's October 2026 Release of AI-Generated Mathematics},pdfauthor={Ying Xie}}

\newcommand{\F}{\mathbb{F}}
\newcommand{\T}{\mathsf{T}}
\providecommand{\lt}{<}
\providecommand{\gt}{>}

\newcommand{\badges}[1]{\par\noindent{\small\sffamily\color{black!70}#1}\par\smallskip}
\newcommand{\badgesep}{\textbar}

\setlist{itemsep=2pt,topsep=3pt}
\title{What the Machine Invented:\\[2pt]
\large Mining New Mathematical Instruments from OpenAI's October 2026 Release of AI-Generated Mathematics}
\author{Ying Xie\\Kennesaw State University\\\texttt{yxie2@kennesaw.edu}}
\date{October 2026}

\makeatletter\let\ps@plain\ps@empty\makeatother
\begin{document}
\maketitle

\begin{abstract}
On 6 October 2026 OpenAI published an initial 722 AI-written mathematical manuscripts, organized into 372 problem families, in the public repository \href{https://github.com/openai/math}{\texttt{github.com/openai/math}}~\cite{src:repo}. The release includes claimed solutions to many prominent open problems. We approach this collection as a treasure mine of mathematical ideas. Our aim is to extract new lemmas, constructions, inequalities and proof mechanisms introduced in the source arguments, and to distinguish them from established tools and new combinations of those tools. We call these \emph{mathematical instruments}. Extraction comes first; their validity, novelty and possible reuse then become subjects for further study. The resulting descriptions also help readers understand the selected papers without having to reconstruct every argument at the outset.
We examine 38 families in two groups of 17 and 21, spanning number theory, geometry, analysis, combinatorics, complexity, algebra, operator algebras, topology and mathematical physics. Each entry records the claim, extracted instruments, antecedents and caveats. Eleven cross-cutting patterns connect the entries, and a scoreboard records formalization status and dependencies. Novelty assessments are provisional and sometimes mixed. Selected source comparisons correct mathematical statements, hypotheses and counts; they do not independently establish priority or certify the underlying proofs. No Lean development was compiled or kernel-checked for this survey.
We also briefly describe a subsequent application of the signed restriction lemma from family~112 to circuit lower bounds at every fixed depth.
\end{abstract}

\tableofcontents
\bigskip

\section{Introduction}\label{sec:intro}

OpenAI's release brought together hundreds of manuscripts addressing longstanding mathematical problems: the Unique Games Conjecture, the exponent of matrix multiplication, Hadwiger's conjecture, the chromatic number of the plane, Kaplansky's conjectures, the amenability of Thompson's group~$F$, the isomorphism problem for free group factors, and many others. Several are accompanied by Lean formalizations; most are not. None had been peer-reviewed at the time of writing.

The motivation for this project is to treat the release as a source of mathematical ideas. A proof is more than a certificate that its theorem is true: it is a store of methods. When a long proof of a famous conjecture is correct, the devices that make it work can have value beyond the theorem; even when an argument fails, a lemma or construction within it may remain valid. The purpose of this project is to mine the release for those devices. We call them \emph{instruments}: new objects, structures, inequalities or proof architectures that a proof introduces in order to reach its conclusion, together with new settings in which established tools are applied.

The project has two goals.
\begin{enumerate}[label=(\roman*)]
\item \textbf{Extracting new instruments.} To isolate and describe the new lemmas, constructions and proof mechanisms generated in these arguments. Once their statements, proofs and relation to earlier work are understood, they can contribute to the shared mathematical knowledge base.
\item \textbf{Making the source arguments accessible.} To give readers a consistent account of what each selected paper claims, which instruments it introduces, which established results it uses, and which steps require closer examination.
\end{enumerate}

\paragraph{A case study in reuse.} After extraction, we investigate whether particular instruments can be reused. A separate application manuscript studies the signed restriction lemma identified in family~\hyperref[fam:112]{112}. The lemma controls correlations through signed combinations of narrow CNFs with bounded expected coefficient mass. Combining it with adaptive averaging, classical depth reduction and quantitative sparsification yields, for the same polynomial-time language as in the source, the gate-count lower bound
\[
 S_D(f_n)\ge 2^{\Omega_D\left((n\log n/\log\log n)^{1/(D-1)}\right)}
 \qquad\text{for every fixed }D\ge3
\]
at all sufficiently large input lengths~$n$. Here $S_D$ counts AND/OR gates in exact circuits of depth at most~$D$, with arbitrary finite fan-in and fan-out and input negations. The source supplies the lemma and the language; the application extends the source's stated depth-three result. The proof will be presented in the separate application paper. This case illustrates the progression from extracting an instrument to investigating its use and obtaining a further theorem.

\paragraph{What this paper contains.} We report an extraction and preliminary novelty assessment for 38 families, selected for the prominence of their problems. For each family, the initial extraction examined the source argument to identify candidate new instruments and distinguish them from combinations of known tools. The paper contains:
\begin{itemize}
\item a \emph{map} of the 38 families (Tables~\ref{tab:map1} and~\ref{tab:map2}), giving for each the result, its Lean status, its central extracted instrument and a provisional novelty assessment;
\item \emph{eleven cross-cutting patterns} (Section~\ref{sec:patterns}) describing recurring approaches in the source proofs, from ``existence instead of construction'' to ``established methods in new settings'' and ``rediscovery without citation'';
\item \emph{family-by-family cards} (Sections~\ref{sec:cards} and~\ref{sec:cards2}) that state each claim, describe the selected candidate instruments with section and lemma references, list what is reused, and record caveats;
\item a \emph{verification and dependency scoreboard} (Section~\ref{sec:scoreboard}) recording formalization status, statement faithfulness, dependence on other OpenAI manuscripts, and uncited nearest ancestors;
\item an explicit account of the \emph{limits of the extraction and mathematical review} (Section~\ref{sec:limits}).
\end{itemize}

\paragraph{Status of the results.} Throughout, ``the paper proves'' and similar phrases describe what a manuscript \emph{claims} and how its argument is organized. The manuscripts have not been peer-reviewed. Formalization statuses are those recorded during the initial extraction as of 11 October 2026, with the source-status corrections noted below. No Lean project was compiled during either stage. Verdicts of novelty mean ``new relative to the citable literature examined during the initial extraction'', not established priority.

\section{Background and motivation}\label{sec:background}

\subsection{The challenge of reading a large release}
A release of this size poses several practical difficulties for readers:
\begin{itemize}
\item \textbf{Volume.} Some manuscripts run to 170--200 pages, and several families contain closely related papers. Reading even a selected part of the release requires substantial time.
\item \textbf{Technical density.} Many arguments involve nested limits, numerous parameters, and compactness or ultrafilter constructions. Identifying the central mechanism can require tracing dependencies across many sections.
\item \textbf{No peer review.} The manuscripts have not passed through the ordinary refereeing process, and the repository's own formalization catalogue marks review status ``unchecked''.
\item \textbf{Need for accessible descriptions.} A short account of an instrument can help readers decide where to begin in a long source argument, while leaving the correctness of that argument to further examination.
\end{itemize}
These difficulties make structured descriptions useful: they provide starting points for reading, comparison and verification.

\subsection{Extracting the instruments within the proofs}
The project begins with the mathematical content inside each selected argument. We identify the objects, lemmas and proof mechanisms it introduces, trace their roles in the proof, and compare them with earlier work. The family entries present the results of that extraction in a common format. Whether a particular instrument can be used elsewhere is a subsequent question, illustrated by the application described in the introduction.

\subsection{Why this matters}
We see several ways in which a structured extraction of this kind can be useful beyond the immediate question of which claims are correct. Some benefits concern access to the source arguments; others, including further mathematical applications, depend on the validity of the instruments concerned.
\begin{enumerate}[label=(\arabic*)]
\item \textbf{A catalogue of extracted instruments.} Instruments are described in a form detached from their headline theorems: for example, a ``shared-leg entropy inequality'' for points of Strassen's asymptotic spectrum (\hyperref[fam:107]{107}), a ``latent alphabet gadget'' with noise that rarely changes the honest decoder's output (\hyperref[fam:102]{102}), ``triangular polynomial cells'' with descending precision budgets (\hyperref[fam:159]{159}), or a ``coupled purity bootstrap'' as a spectral-gap criterion (\hyperref[fam:268]{268}). Their statements and proofs can be studied separately from the headline claims before possible applications are investigated.
\item \textbf{Cross-field transfer of mechanisms.} One recurring pattern is ``established methods in new settings'' (Section~\ref{sec:patterns}): a published theorem applied in a setting for which the initial comparison found no precedent, such as the Benyamini--Sternfeld fixed-point-free Lipschitz map against amenability, Braverman's theorem on bounded-depth circuits for Chowla, or Quot-scheme localization for Bloch's conjecture. Describing these transfers gives readers concrete examples to investigate in other settings.
\item \textbf{Directing verification effort.} The cards and scoreboard identify steps on which the conclusions depend, dependencies on other unverified manuscripts, and gaps in the formalization coverage inspected. These details help a reader choose which statements or arguments to check first.
\item \textbf{Reducing the cost of initial reading.} A short entry introduces a manuscript that may extend to two hundred pages. It points to the proposed new instruments and distinguishes them from established ingredients, giving a specialist a route into the full proof.
\item \textbf{A method for examining further releases.} The same extraction method can be applied to other large collections. A fixed protocol (one AI reading per family, usually up to three candidate devices traced through the proof, a literature search per device, a four-grade verdict, and an explicit list of caveats and dependencies) is a template that the community can apply, extend, and audit in turn.
\item \textbf{Scaffolding for formalization.} By isolating the exact statement, the load-bearing lemmas, and the external and internal dependencies of each result, the extraction indicates where formalization effort should be directed, and which existing formal developments still require independent compilation and kernel re-checking.
\end{enumerate}

\subsection{What the extraction does not do}
The extraction is not a verification of the source proofs; selected statements and identities were checked. Its novelty assessments are limited by the literature examined, and fresh searches were not completed for every family in the second group. The selection covers 38 of 372 families and favors prominent problems. The instruments identified here are devices within proofs; the extraction did not identify a new foundational theory comparable in scope to motives or perfectoid spaces. Section~\ref{sec:limits} gives the full limits of the extraction and mathematical review.

\section{The extraction at a glance}\label{sec:glance}
\noindent\textit{OpenAI math release \textperiodcentered{} 6 Oct 2026 \textperiodcentered{} \textbf{extraction from 38 selected families}.}\par\medskip

The selected set contains thirty-eight family entries, in groups of seventeen and twenty-one. The initial extraction asks which instruments appear new and which combine known tools. The tables retain its qualitative verdicts, with corrected theorem statements and explicit limits on verification.

\begin{table}[ht]\centering\small
\begin{tabular}{@{}>{\raggedright\arraybackslash}p{0.2\textwidth}>{\raggedright\arraybackslash}p{0.15\textwidth}>{\raggedright\arraybackslash}p{0.57\textwidth}@{}}
\toprule
\textbf{Quantity} & \textbf{Value} & \textbf{Note}\\ \midrule
Families examined & \textbf{38} & selected for the prominence of the problem, across number theory, geometry, analysis, combinatorics, complexity, algebra, operator algebras, topology and mathematical physics. Group 1: 17 families. Group 2: 21 families \\ \addlinespace
Novelty classifications & qualitative & Mixed labels are retained at the family and instrument level. No disjoint numerical breakdown is asserted, and novelty judgments have not been independently validated by the subsequent mathematical review \\ \addlinespace
Original headline-formalization labels & \textbf{19 of 38} & Includes partial inspection for 003 and 087. Five other original entries were companion-only, counting 159. These are labels from the initial extraction; source-status corrections and the counting rule appear in Section~\ref{sec:scoreboard}. No Lean build was run \\ \addlinespace
Recurring pattern & \textbf{Established methods in new settings} & a 1983 theorem on maps without approximate fixed points is applied to amenability; a 2009 circuit-complexity theorem to Chowla correlations; Quot-scheme localization to Bloch; a 2019 conformal map to Mahler \\ \addlinespace
\bottomrule
\end{tabular}
\caption{Summary figures of the extraction (as of 11 Oct 2026).}\label{tab:glance}
\end{table}
\section{\texorpdfstring{Method, and what ``new'' means here}{Method, and what new means here}}\label{sec:method}

The work proceeded in two stages. The \emph{initial extraction} identified candidate mathematical instruments, compared them with the literature, and assigned preliminary novelty assessments. The \emph{subsequent mathematical review} checked the reported theorem statements and selected supporting arguments and revised the draft. These terms distinguish the two activities throughout the paper. Both groups of families---the first 17 and the second 21---belong to the initial extraction.

During the initial extraction, a separate AI reading was carried out for each family, using the LaTeX source where the repository provides it (most do, under \texttt{build/\allowbreak{}sections/\allowbreak{}}), otherwise the PDF, together with the Lean challenge and solution files and, where released, the abridged reasoning trace. Each reading generally isolated up to three candidate devices, traced each through the proof, then searched the literature for the device and its nearest relatives. Verdicts use four grades.

\begin{description}
\item[New] an object, structure, inequality or proof architecture with no precedent found, even though its parts are classical.
\item[New combination] the ingredients are established, while their combination or application appears new.
\item[Incremental] a refinement or variation of an established technique.
\item[Unclear] the available text and comparisons did not support a definite classification.
\end{description}

The subsequent mathematical review compared theorem statements and selected supporting arguments against a pinned repository snapshot and corrected mathematical statements, hypotheses and counts. Potential extensions, including the proposed logarithm analogue of family 017, are distinguished from consequences actually established here. Literature comparisons and novelty assessments originate in the initial extraction; the subsequent mathematical review did not repeat every search.

Literature-search coverage varied across the two groups. For some entries in the second group, comparison relied on directly retrieved primary sources and the model's existing knowledge without a fresh search. The affected entries identify this limitation. These comparisons are not an exhaustive prior-art search.

\section{\texorpdfstring{The map}{The map}}\label{sec:map}

{\footnotesize\setlength{\tabcolsep}{3pt}
\begin{longtable}{@{}>{\raggedright\arraybackslash}p{0.065\textwidth}>{\raggedright\arraybackslash}p{0.17\textwidth}>{\raggedright\arraybackslash}p{0.10\textwidth}>{\raggedright\arraybackslash}p{0.44\textwidth}>{\raggedright\arraybackslash}p{0.15\textwidth}@{}}
\caption{The map, first group: 17 families.}\label{tab:map1}\\
\toprule
\textbf{Family} & \textbf{Result} & \textbf{Lean} & \textbf{Central extracted instrument} & \textbf{Verdict} \\
\midrule
\endfirsthead
\multicolumn{5}{@{}l}{\emph{(continued)}}\\
\toprule
\textbf{Family} & \textbf{Result} & \textbf{Lean} & \textbf{Central extracted instrument} & \textbf{Verdict} \\
\midrule
\endhead
\bottomrule
\endfoot
\hyperref[fam:003]{003} & $\zeta$ and Dirichlet L zero-free for $\operatorname{Re} s>7/8$ & partial inspection & Twist the Möbius sum by all sextic residue symbols; Poisson fuses $\mu$ with Gauss sums into cubic Gauss sums, whose Kubota--Patterson automorphy moves the thin sextic family to a quadratic one where square-root cancellation on average is known. & combination as architecture, leaning new \\ \addlinespace
\hyperref[fam:004]{004} & Hilbert's tenth problem over $\mathbb{Q}$ undecidable & no & The reduction does not construct a Diophantine definition of $\mathbb{Z}$. A Turing reduction dovetails integer-root enumeration with finite rational-oracle tests. Compactness and elliptic height bounds show that some finite test fails if no integer root exists. Every terminating run makes finitely many queries, with no uniform bound asserted. & new \\ \addlinespace
\hyperref[fam:017]{017} & Irrationality exponent of $\pi$ equals 2 & yes & ``Separated-weight interpolation'': jets at translates on the exponential curve with anisotropic weights chosen independently of the centers, while the sufficiently large degree threshold may depend on them, proved by a Dyson-type curve inequality pushed through nef/ample and Serre vanishing on a weighted blow-up. & combination bordering new \\ \addlinespace
\hyperref[fam:032]{032} & Rational Hodge conjecture for CM abelian varieties & no & ``Fixed-input variation'': freeze the theta-series input on a Picard-type ball quotient, move only the definite-place complex structure, keep a (2,0)-form period nonzero, then rationalize by weak approximation. & combination \\ \addlinespace
\hyperref[fam:074]{074} & Kakeya maximal conjecture in $\mathbb{R}^{3}$; Kakeya sets in $\mathbb{R}^{4}$ have dimension 4 & no & Differentiate the critical exponent: BV gives approximate differentiability of the exponent surface; a separate monotonicity argument yields an ordinary differentiability point whose partial derivatives prescribe extremal configurations ruled out by an entropy and mutual-information bound. In 4D, a variational ``horizon'' built from quadratic local fits replaces polynomial Wolff axioms. & combination leaning new; 4D framework new \\ \addlinespace
\hyperref[fam:087]{087} & Mahler conjecture, symmetric and general & symmetric & Replace each real strip by a complex lens with uniform boundary coordinate (Gross's conformal Skorokhod map); push Lelong-number Jacobian mass to the boundary as $k\to\infty$; the volume product becomes a sum of feasibility probabilities over signed simplices packing inside the polar. & new feasible-simplex identity; rest combination \\ \addlinespace
\hyperref[fam:102]{102} & Unique Games Conjecture & yes & ``Latent alphabet gadget'': a nonlinear shift-equivariant decoder on a big latent space with a noise law that rarely changes the honest decoder but is detected by every linear family with sufficiently large rank on the output subspace, so the degree-two KMS inverse theorem applies unchanged. & new \\ \addlinespace
\hyperref[fam:107]{107} & Matrix multiplication exponent $\omega \le 9/4$ & yes & The proof bounds every point of Strassen's asymptotic spectrum using polynomial multiplication, without giving a matrix-multiplication algorithm: a DFT-over-copies degeneration keeps all blocks sharing a leg and yields a ``shared-leg entropy inequality''; concavity plus tripling force growth $a^{4/3}$. & new as a proof of an $\omega$ bound \\ \addlinespace
\hyperref[fam:157]{157} & Counterexample to Hadwiger's conjecture & yes & Algebraic ``hole relation'': non-edges certified by a symmetric bilinear identity whose terms cancel around any triangle leaving $0=1$, so the complement is triangle-free. Connected matchings are excluded by an entropy-minimizing container argument. & new \\ \addlinespace
\hyperref[fam:158]{158} & The plane is not 5-colorable & yes & ``Transfer of colorability'': a proper k-coloring exists iff a weak measurable one does, via a rigidity theorem that every wild rotation-invariant law on the dual of the algebraic plane is Haar. The only finite graph used is a Moser spindle in the final step. & new \\ \addlinespace
\hyperref[fam:159]{159} & Erd\H{o}s reciprocal-sum conjecture; quasipolynomial Szemerédi & original: companion only; see card & ``Triangular polynomial cells'' with descending precision budgets: a density increment that never resets to a sub-progression; the loss at weight i depends only on precisions at higher weights, turning $\log\log N$ into $\log N$. & new architecture \\ \addlinespace
\hyperref[fam:196]{196} & Counterexample to Kaplansky's zero-divisor conjecture & yes & Projective-plane label sets with pairwise odd intersections make every vertex of the ``simultaneous-step'' product graph odd-degree, so every component is even and the product vanishes mod 2. Collapse is excluded by a girth-conditioned random-matching estimate plus a Lipton--Tarjan reduction of planar pairings. & combination; separator step new \\ \addlinespace
\hyperref[fam:197]{197} & Counterexample to Kaplansky's direct-finiteness conjecture; a non-sofic group as a consequence & yes & An ``incidence criterion'': a hypergraph of degree $\ell+1$ and girth 12 plus vectors with block-sum 0 and total $-I$ yields $ab=1,\ ba\ne1$ after gluing copies of $\F_\ell^2$, averaging, character-twisted HNN packing, and a homomorphism to the field that exposes a rank-deficient Gram matrix. & new criterion; parts combination \\ \addlinespace
\hyperref[fam:248]{248} & Thompson's group F is not amenable & yes & Color each element of a Følner set by a Hilbert-ball vector via the recursion ``color = f(average of child colors)'', with f a Benyamini--Sternfeld Lipschitz map having no approximate fixed point; pair transports force the averages to agree, an impossibility. & new use; leaning new \\ \addlinespace
\hyperref[fam:287]{287} & All free group factors isomorphic; fundamental group $\mathbb{R}_{+}$ & yes & ``Small cocycle with a prescribed word value'': a Fox-derivative cocycle $\ell^{2}$-tiny on each generator that accumulates a full unitary factor on a long word with free prefixes, because spectral cutoffs of the prefix-sum operator give approximate inverses in the required vector norm (the limiting free Poisson law has no atom at 0). A trace-preserving flow then absorbs a generator. & combination \\ \addlinespace
\hyperref[fam:304]{304} & Hilbert--Smith conjecture in every dimension & no & Witt-group ``signature'' of a sphere's orientation form pushed through the wild orbit space into a test sphere, in ``compact-source test categories'' of sheaves built from arbitrary continuous maps; denominators bounded by a lattice fixed before any action data; a faithful $\mathbb{Z}_p$ action splits it into $p^k$ equal integral pieces. & new \\ \addlinespace
\hyperref[fam:362]{362} & 3D relativistic Vlasov--Maxwell, large data, global smooth & yes & ``Signed momentum increment'' bootstrap: bound the net impulse on a particle over a short interval, not the force pointwise, via a signed Liénard--Wiechert identity along source worldlines; each momentum doubling costs time $\ge c/\log P$, and the harmonic sum diverges. & combination; direction-change counting new \\ \addlinespace
\end{longtable}}

\section{\texorpdfstring{The map, second group}{The map, second group}}\label{sec:map2}

{\footnotesize\setlength{\tabcolsep}{3pt}
\begin{longtable}{@{}>{\raggedright\arraybackslash}p{0.065\textwidth}>{\raggedright\arraybackslash}p{0.17\textwidth}>{\raggedright\arraybackslash}p{0.10\textwidth}>{\raggedright\arraybackslash}p{0.44\textwidth}>{\raggedright\arraybackslash}p{0.15\textwidth}@{}}
\caption{The map, second group: 21 families.}\label{tab:map2}\\
\toprule
\textbf{Family} & \textbf{Result} & \textbf{Lean} & \textbf{Central extracted instrument} & \textbf{Verdict} \\
\midrule
\endfirsthead
\multicolumn{5}{@{}l}{\emph{(continued)}}\\
\toprule
\textbf{Family} & \textbf{Result} & \textbf{Lean} & \textbf{Central extracted instrument} & \textbf{Verdict} \\
\midrule
\endhead
\bottomrule
\endfoot
\hyperref[fam:002]{002} & Full BSD formula in Selmer corank 0 and 1 & no & ``Tame variables'': ultrafilter limits of cyclic characters at auxiliary primes used as extra Iwasawa variables to force integrality with no residual-image hypothesis; a ``theta implication'' comparing Artin lengths from a holomorphic theta family on U(3,1) with a nilpotent parameter replaces Eisenstein congruences. & new construction; rest combination \\ \addlinespace
\hyperref[fam:007]{007} & Ordinary two-point Chowla; corrected binary Elliott & original: statements only; see card & ``Terminal-scale divisor bins'': dilations confined so narrowly that the scale is unchanged, removing Tao's logarithmic averaging; integers certified to behave like independent residues by Braverman's theorem that polylog independence fools $\mathrm{AC}^{0}$. & combination bordering new \\ \addlinespace
\hyperref[fam:010]{010} & Fontaine--Mazur modularity at $p=2$, all residual types & no & ``Marked lines'': after a huge abelian 2-extension, record an ordinary line at only a vanishing fraction of dyadic places, enough to force the residual fiber connected while leaving the propagation step unrestricted at 2. & combination \\ \addlinespace
\hyperref[fam:040]{040} & Bloch's conjecture for all surfaces with $p_g=0$ & no & Apply virtual enumerative geometry to a Chow problem: a torus-localized Quot-scheme count with two copies of the surface left unintegrated vanishes for degree reasons, so fixed loci give a diagonal decomposition with a universal coefficient; for fake projective planes the coefficient is read on a Cartwright--Steger surface. & new strategy \\ \addlinespace
\hyperref[fam:103]{103} & $\mathrm{L}=\mathrm{RL}=\mathrm{BPL}$ & yes & ``Exact correction and copy hierarchy'': no pseudorandom generator; refine the configuration graph algebraically so every transition weight is small, estimate the reward with a sampler driven by one O(log n)-bit environment, and control the evaluator's recursive space use with a catalytic key-equality gadget. & new \\ \addlinespace
\hyperref[fam:112]{112} & Explicit depth-3 circuit lower bound $2^{\omega(\sqrt{n})}$ & yes & ``Restriction bound independent of clause count'': a restricted width-k CNF written as a signed combination of width-b CNFs with bounded expected mass, with a suitable testing width for any fixed $pk$; the standard switching bound $(5pk)^s$ gives decay in $s$ when $pk<1/5$. & new \\ \addlinespace
\hyperref[fam:143]{143} & Hilbert's 16th: $H(n)<\infty$ for every degree & no (companion only) & Parametric Lion--Rolin preparation of the slope field, coefficients as parameters, in place of blow-ups, giving a finite alphabet of passages uniform over degeneracies; a Khovanskii-type projection principle turns non-uniform finiteness into a uniform bound; Ilyashenko's cochains re-engineered into multi-scale ``packets''. & combination \\ \addlinespace
\hyperref[fam:156]{156} & Borsuk's conjecture fails in dimension 9 & yes & The set of rank-one projectors of $\mathbb{R}^{4}$ (a copy of $\mathbb{R}P^3$ in $\mathbb{R}^{9}$) needs at least 11 pieces: an ``odd spectral extension'' turns any 10-cover into an odd map between 9-spheres, and degree analysis at the tight Borsuk--Ulam bound forces a 10-vertex $\mathbb{R}P^3$-like complex that cannot exist. & combination leaning new \\ \addlinespace
\hyperref[fam:193]{193} & Serre positivity in ramified mixed characteristic & no & ``Rank quotient'': sequences of finite-length modules along the p-power-root tower modulo negligible normalized rank, packaged as perfect complexes over a connective algebra so K-theoretic descent acts on asymptotic perfectoid data; ordinary and perfectoid Euler characteristics compared only after Galois norms. & new device \\ \addlinespace
\hyperref[fam:221]{221} & Mézard--Parisi formula for diluted spin glasses & yes (not in catalogue) & ``Depth shift'': couple to a finite auxiliary hierarchy by marked Poisson perturbations and shift every target branching depth by one level so divisions by small coefficients stay bounded, proving the factorization Panchenko left open. & new \\ \addlinespace
\hyperref[fam:263]{263} & Ionization conjecture, full Schrödinger theory & companions only & Condition a many-body ground state on noisy, unlabeled coarse-grained observations of electron positions; the energy cost is logarithmic in the event probability and free of particle number, so local failures of Thomas--Fermi yield low-energy states that contradict semiclassics. & new \\ \addlinespace
\hyperref[fam:265]{265} & 2D area law from a global gap, frustration allowed & no & ``Optimized nested positive filters'' whose only cost is excitation energy, and ``replica metric transport'' encoding conditional mutual information as the derivative of a log-norm comparable to energy, with a displayed information term linear in the transport parameter and additional energy cost quadratic in that parameter and linear in the split probability. & new \\ \addlinespace
\hyperref[fam:266]{266} & Exactly three MUBs in dimension six & headline no; $N(6) \le 5$ yes & Represent every order-6 complex Hadamard matrix by nine seed phases plus a U(2) polar factor, giving a finite cover of the unclassified continuum; the headline claim uses a binary64 computation with hand-proved margins. Companion: ``injective moments'' calculus establishing Fourier vanishing outside Tao's $S_6$ equivalence class. & combination \\ \addlinespace
\hyperref[fam:268]{268} & Haldane gap for the spin-1 Heisenberg ring & no & ``Coupled purity bootstrap'': one partition function written as Gibbs weights and as powers of spatial transfer eigenvalues; squaring one doubles inverse temperature, squaring the other doubles length, and exact identities couple them. The criterion does not require frustration-freeness. & new criterion \\ \addlinespace
\hyperref[fam:285]{285} & Reduced Baum--Connes injectivity fails (torsion-free) & no & ``Heisenberg voltage tower'': relators are nilpotent-group covers of an expander tower so each layer carries $\ell^2(H)$ and incidence operators realize Bott-class symbols; an Eilenberg swindle makes the embedded torus class vanish. & combination \\ \addlinespace
\hyperref[fam:293]{293} & No hyperinvariant subspace; transitive algebra problem & yes & ``Backward intertwiner over an ergodic base'': unicellular weighted shifts glued over the 2-adic odometer; multipliers force a tail index per fiber, a commuting operator lowers it along orbits, and a finite integer-valued measurable label cannot strictly decrease almost everywhere under a measure-preserving map. & new \\ \addlinespace
\hyperref[fam:297]{297} & Naimark's problem in ZFC & yes (not in catalogue) & ``Selective bridge'': a reduced HNN adjunction along peaking projection sequences merges exactly two pure-state classes while every other pure state extends uniquely, so each step depends on countably many earlier ones and diamond is replaced by cardinal bookkeeping at length $\mathfrak{c}^{+}$. & combination \\ \addlinespace
\hyperref[fam:309]{309} & Kervaire invariant at $p=3$: survivors exactly $j\in\{0,2,3\}$ & no & ``Finite Euler bound'': a coordinate-change coefficient as a nilpotence witness imposes a deadline on the Tate spectral sequence, replacing the Gap and Periodicity theorems of Hill--Hopkins--Ravenel; the Belmont--Ray coefficient model is proved along the way. & combination \\ \addlinespace
\hyperref[fam:325]{325} & Complete (matrix-valued) Crouzeix conjecture & yes & ``Ordered density comparison'': compress the Crouzeix--Palencia positive kernel by the top singular pair into a $2\times2$ block of boundary matrix functions and measure in Hilbert--Schmidt norm, so Parseval over Faber coefficients replaces the triangle inequality that gave $1+\sqrt{2}$. & combination leaning new \\ \addlinespace
\hyperref[fam:340]{340} & Counterexample to the nearby Lagrangian conjecture & no & ``Momentum separation'': a generating family with an exponential scale over the parameter sphere upgrades the immersed Legendrians of the 2026 Álvarez-Gavela--Igusa--Sullivan program to an embedded exact Lagrangian; the obstruction is a 2-torsion class in Waldhausen's stable tube space. & combination \\ \addlinespace
\hyperref[fam:376]{376} & Universal computation in forced Navier--Stokes flows & lemmas only & Defining the force as the residual enforces Navier--Stokes for a chosen smooth divergence-free flow; the required mean, regularity and time-dependence of that force still need proof; the construction is Moore's generalized shifts plus Bennett reversibility plus Hamiltonian box-moving. The extracted construction uses scale-nested memories under rapidly vanishing forcing. & incremental \\ \addlinespace
\end{longtable}}

\section{\texorpdfstring{Cross-cutting patterns}{Cross-cutting patterns}}\label{sec:patterns}

Several approaches recur across the selected proofs. The following patterns describe how the arguments establish existence, control parameters, transfer results between settings, and introduce auxiliary structures. They provide a way to compare instruments arising in different subjects.

\subsection{\texorpdfstring{Existence instead of construction}{Existence instead of construction}}

Several arguments establish existence through compactness, spectral duality, random sampling, or a terminating search whose output has not been computed. The resulting existence statement need not supply an object of practical size.

The proof of $\omega \le 9/4$ bounds the asymptotic spectrum without providing a competitive finite matrix size. The Hadwiger argument selects m through an ultrafilter, and the Hilbert 10 reduction uses a nonstandard root. The zero-divisor and direct-finiteness constructions use random choices or searches with parameters such as $q=2^{217}$ and $n\gtrsim e^{300000}$. In the second group, the Hilbert 16 argument gives an ineffective B(d), the nearby-Lagrangian example has no explicit dimension, and the Borsuk construction supplies a continuous set without extracting a finite subset. These examples differ in effectivity, but each leaves some construction data unspecified or unevaluated.

\subsection{\texorpdfstring{Estimates become extremal problems}{Estimates become extremal problems}}

A proposed inequality can be studied through the regularity of its optimal-exponent function. Assuming that the desired bound fails, differentiation of this function imposes conditions on an extremal configuration. The proof then rules out a configuration satisfying those conditions.

In the three-dimensional Kakeya argument, coordinatewise monotonicity gives local BV regularity and approximate differentiability almost everywhere. An additional monotonicity comparison yields ordinary total differentiability, and the partial derivatives determine the relevant plank shapes. The four-dimensional argument uses a ``horizon'' defined by an infimum over quadratic fits. In the Erd\H{o}s argument, a potential controls the loss at each weight using only the higher weights.

\subsection{\texorpdfstring{Entropy and information as proof devices}{Entropy and information as proof devices}}

Mutual information, relative-entropy projections, and entropy inequalities organize estimates in several proofs. The information quantities measure how much dependence or conditioning is available, and incompatible upper and lower bounds produce a contradiction.

The Kakeya proof uses a mutual-information bound to apply a projection theorem to dependent events. The Hadwiger construction combines Kleitman--Winston recording with a Csiszár I-projection as its progress measure and a Ford--Fulkerson cut as its stopping certificate. The matrix-multiplication argument for $\omega$ uses an entropy inequality for a shared tensor leg, while the zero-divisor construction bounds repeated traversals of a random graph by entropy. In the second group, conditioning in the ionization argument costs $\log^5(e/p)$ times the inverse squared noise scale in energy. The area-law proof expresses conditional mutual information through a log-norm derivative, and the Haldane-gap proof uses a Rényi-2 purity bootstrap.

\subsection{\texorpdfstring{Parity and characters over finite fields}{Parity and characters over finite fields}}

Four constructions use cancellation in finite fields, either through a scalar identity or through vanishing character sums.

The direct-finiteness construction requires $\ell+1=0$ in K. The zero-divisor construction uses projective lines whose pairwise intersections have odd size. The UGC gadget is defined over $\F_2^s$ and contains a quadratic block over $\F_{2^d}$. The proof of $\omega \le 9/4$ applies a DFT across $L=5M$ copies to select exactly the matching labels.

\subsection{\texorpdfstring{Rigidity theorems that erase regularity hypotheses}{Rigidity theorems that erase regularity hypotheses}}

Some arguments replace a regularity assumption with a transfer theorem. The theorem connects the unrestricted object to a measurable, constructible, or otherwise controlled formulation, so that the remaining proof can use that formulation.

For plane coloring, the transfer states that a proper coloring exists exactly when a weak measurable coloring exists. The Hilbert--Smith argument uses sheaf categories generated by arbitrary continuous maps from polyhedra while retaining invertible Verdier duality. The Hilbert 10 argument uses an oracle reduction without constructing a definable copy of $\mathbb{Z}$.

\subsection{\texorpdfstring{Established methods in new settings}{Established methods in new settings}}

A recurring contribution is the application of an established theorem in a different mathematical setting. The relevant novelty question concerns that application and its surrounding argument; the imported theorem retains its original attribution.

The examples include the Benyamini--Sternfeld theorem of 1983 on Lipschitz maps without approximate fixed points in an amenability argument; Gross\textquoteright{}s 2019 conformal Skorokhod map in the Mahler proof; and Kubota--Patterson cubic theta theory from 1977 in a zero-free-region argument. Other transfers use Fox calculus and free Poisson laws for free group factors, Lelong numbers for volume products, and Lipton--Tarjan separators for diagrams of van Kampen type. The second group applies Braverman\textquoteright{}s 2009 $\mathrm{AC}^{0}$ theorem to Chowla correlations, Marian--Oprea--Pandharipande localization to Bloch\textquoteright{}s conjecture, Lion--Rolin preparation to Hilbert\textquoteright{}s 16th problem, Ueda\textquoteright{}s reduced HNN extensions to Naimark\textquoteright{}s problem, Suzuki\textquoteright{}s quantum transfer matrix to the Haldane gap, and Kubo--Ando means to the area law.

\subsection{\texorpdfstring{Rediscovery without citation}{Rediscovery without citation}}

The initial comparison identified several constructions closely related to earlier work that was not cited in the manuscript. In some cases, the corresponding reasoning summary does name the earlier source. These comparisons matter when distinguishing a new instrument from a reconstruction or adaptation.

The Vlasov--Maxwell manuscript uses a time-integrated momentum increment related to Cheverry--Ibrahim's work, which is cited six times in the reasoning summary but absent from the paper. The interpolation argument for $\pi$ has antecedents in Masser and Fischler--Nakamaye; the spectral argument for $\omega$ follows a template in Alman--Li's 2026 ``Strassen calculus''; and the Kakeya comparison identifies Rios--Sawyer as an uncited antecedent. In the second group, the M\'ezard--Parisi reasoning summary names Mourrat--Panchenko and Barbier--Panchenko, while the paper omits them. The Serre-positivity argument uses an orbit lemma of Suslin--Voevodsky without citation.

\subsection{Consequences and proposed extensions}

Some entries describe consequences of the stated results or suggest extensions beyond them. Each requires its own hypotheses and argument. The cards distinguish consequences presented as following from the source from adaptations that still need a proof.

The centers $(r^j,jp_i/q_i)$ suggest a possible analogue of the $\pi$ argument for $r\in\mathbb Q_{>0}\setminus\{1\}$. Establishing $\mu(\log r)=2$ would require a separate interpolation, height and determinant argument; it does not follow verbatim from the theorem surveyed here. The claimed free-group-factor isomorphism would refute the Connes--Shlyakhtenko $L^2$-Betti conjecture for $L(\F_n)$, a consequence not discussed in the paper. In the second group, the hyperinvariant-subspace counterexample would also refute Kadison's transitive-algebra problem, which the catalogue summary omits; the Mézard--Parisi lower bound is stated for bounded interactions with permutation-invariant law, odd arity included.

\subsection{\texorpdfstring{Probability imported into deterministic problems}{Probability imported into deterministic problems}}

Several proofs introduce a random or ergodic auxiliary system into a deterministic problem. Almost-sure statements or measure-zero obstructions in that system then support the deterministic conclusion.

In the hyperinvariant-subspace construction, operators are assembled over the 2-adic odometer, where a finite integer-valued measurable label cannot strictly decrease almost everywhere along the measure-preserving dynamics. The ionization argument smooths and conditions electron positions. Marked Poisson perturbations provide the hierarchy in the Mézard--Parisi proof, while the Baum--Connes construction uses a local-lemma probability law along an infinite tower. The plane-coloring argument uses Haar rigidity for the unrestricted laws considered there.

\subsection{\texorpdfstring{Auxiliary structure manufactured by the proof}{Auxiliary structure manufactured by the proof}}

An argument may introduce auxiliary structure with properties it can control, and then prove that the original object is constrained by that structure.

The UGC construction introduces a latent alphabet hidden from the honest player. The proof of $\mathrm{L}=\mathrm{BPL}$ uses one shared O(log n)-bit environment for all samplers. In the Mézard--Parisi argument, a finite auxiliary hierarchy is refined only after $N\to\infty$. The ionization proof uses a decreasing chain of observation $\sigma$-fields. The BSD argument introduces tame Iwasawa variables and specializes them to zero at the end.

\subsection{\texorpdfstring{Label versus content}{Label versus content}}

A title may suppress restrictions that are essential to the mathematical statement. The extraction therefore records the precise model, hypotheses, and conclusion, including computational assumptions where they occur.

``Universal computation in Navier--Stokes'' uses volume-preserving flows and a force defined by the residual; the prescribed forcing properties, including eventual stationarity, remain mathematical constraints to verify. ``Baum--Connes counterexample'' is reduced, coefficient-free rational injectivity for a finitely generated group; finite presentation is not determined, and the claim does not refute the classical Novikov conjecture. ``Three MUBs in dimension six'' rests on a binary64 computation conditional on compiler semantics. ``Ionization conjecture'' is the strict-binding formulation, with threshold bound states not excluded. ``Beyond the square-root exponent'' is $2^{\omega(\sqrt{n})}$ at inverse-tower rate, not $2^{n^{1/2+\varepsilon}}$ for a fixed $\varepsilon>0$.

\section{\texorpdfstring{Family by family}{Family by family}}\label{sec:cards}

Each card states the claim, describes the extracted instruments and their roles, lists established ingredients, and records caveats. Section and lemma numbers refer to the source papers. The compact status lines report supplied formalizations or the initial extraction's assessment; they do not report a Lean build performed for this survey.

\subsection{Number theory and arithmetic geometry}

\subsubsection{\texorpdfstring{003 \textperiodcentered{} Zero-free half-plane $\operatorname{Re} s>7/8$}{003 · Zero-free half-plane Re s > 7/8}}\label{fam:003}
\badges{Lean, agent-generated, unreviewed \badgesep{} architecture: combination, leaning new}

Every Dirichlet L-function and every finite-order Hecke L-function over $\mathbb{Q}(\sqrt{-3})$ has no zeros in $\operatorname{Re} s>7/8$. The main paper has 199 pages; a human-edited companion reaches 11/12. Corollaries include Vinogradov's least-nonresidue conjecture.

\paragraph{Extracted instruments.}

\begin{itemize}
\item \textbf{Embedding the Möbius sum in a sextic family.} Set $A_u(D)=\sum\mu(n)\nu(n)\chi_n(u)W(Nn/D)$ with $\chi_n$ the sextic symbol. Rows $u=p^6$ are near-copies of the target, so a mean-square bound over $u$ controls $A_1$. Poisson in $u$ produces sextic Gauss sums; the Hasse--Heath-Brown identity absorbs $\mu$ into cubic Gauss sums, which are Fourier coefficients of the Kubota--Patterson cubic theta function; its transformation law converts the sextic twist into a quadratic one, where Goldmakher--Louvel's large sieve has no $(ML)^{2/3}$ loss. The general sextic sieve would only give $D^{7/3}$. The initial search found no earlier use of cubic metaplectic automorphy to bound M\"obius sums.
\item \textbf{Cube completion and transfer.} Keeping the full index $nb^3$ allows automorphy to be applied. The argument then removes cubes by M\"obius inversion and iterates with contracted row ranges. The initial extraction classifies this as an incremental development of Heath-Brown 1995 and Dunn--Radziwi\l{}\l{}.
\item \textbf{Continuation from a common signal, saturated-witness zero detector, compensated probe} (Part II, the push from 11/12 to 7/8). This introduces additional accounting into Littlewood's $\mu\Rightarrow{}$zero-free reduction and the Montgomery / Maynard--Pratt detector; the initial assessment is incremental.
\end{itemize}

\paragraph{Reused.}

Kubota--Patterson theta and Dunn--Radziwi\l{}\l{} cusp expansions; Goldmakher--Louvel and Baier--Bansal large sieves; Heath-Brown's recursion; Thorner--Zaman prime counting; Hecke functional equation.

\paragraph{Caveats.}

The needed bound is comparable to GRH on average over a family for which the general cubic sieve is sharp. The twist bookkeeping in the theta reflection is therefore a central point to check, as is the dependence of constants on the character.

\subsubsection{\texorpdfstring{004 \textperiodcentered{} Hilbert's tenth problem over $\mathbb{Q}$}{004 · Hilbert's tenth problem over ℚ}}\label{fam:004}
\badges{not formalized \badgesep{} new architecture \badgesep{} depends on two other OpenAI papers}

No algorithm decides whether an integer polynomial has a rational zero; the problem has Turing degree $0'$, already for degree $\le 4$. The paper explicitly gives no existential definition of $\mathbb{Z}$ in $\mathbb{Q}$ and no Diophantine model, so this route does not conflict with Mazur's conjecture.

\paragraph{Extracted instruments.}

\begin{itemize}
\item \textbf{Finite rational tests with compactness.} The construction starts with a recursive theory in which every ring element $a$ carries an elliptic point $T_a\in mE(F_R)$, $F=\mathbb Q(\sqrt2)$, on a fixed rank-one curve. Adjoin $f(a)=0$, Skolemize, and feed finite initial segments of the ground constraints to the $\mathbb{Q}$-oracle. Dovetail these tests with enumeration of integer tuples: a root ends the enumeration, while failure of a test ends the no-root branch. Each terminating run uses finitely many oracle queries, without a uniform bound. To prove termination, suppose all tests pass; compactness embeds a model in the rational sort of an elementary extension; ``elliptic indices'' $\eta:R\to{}^*\mathbb Z$ are additive and injective; local tests at ``contact primes'' force $\tilde a\equiv\eta(a)$ to high $q$-adic precision; a height bound $h(\tilde A)\ll_f B$ against $h(x(\eta(a_j)P))\approx2\hat h(P)\eta(a_j)^2$ gives $B^2\ll B$, so the indices are ordinary integers and the root is integral. This replaces the Diophantine-definition approach (Poonen, Koenigsmann, Mazur--Rubin, Koymans--Pagano). The initial comparison found no matching architecture.
\item \textbf{Five-point height estimate.} For the source's fixed five-point set $\mathcal C$, fixed exceptional prime set $S$, and constants $H,c>0$, the bound is $h(s)\le H\cdot M_S(s)^c$ for $s\in\mathbf P^1(\mathbb Q)\setminus\mathcal C$. Here $h$ is logarithmic Weil height, and $M_S$ is the squarefree product of primes outside $S$ with odd contact order at one of the five fixed branch values of the Shimura curve $X=\mathrm{Sh}(B_{210})$, genus 5, modulo its Atkin--Lehner group. Even contact gives an unramified fiber; a sign-obstruction splitting yields an odd 2-adic representation; modularity, Javanpeykar's bound, Gaudron--Rémond isogeny estimates and Pazuki's height comparison complete the argument. The statement concerns this specified five-point configuration; no analogue for other branch sets is asserted here. The initial extraction classifies the statement as new and the proof as a combination of established methods.
\item \textbf{Pole-parity existential formula.} The construction gives $\Phi$ with $\mathbb{Z}\subseteq\Phi(\mathbb{Q})\subseteq\{b:{}$ negative valuations are even off a finite set$\}$, via 2-descent on two curves with Monsky matrices and exact Selmer rank 3. This replaces the Demeyer--Van Geel and quaternion predicates and is classified as a new combination.
\end{itemize}

\paragraph{Reused.}

DPRM; Poonen 2002 index reading; Cassels--Tate; Deligne canonical models; Darmon--Granville; Green--Tao--Ziegler with Selberg sieve; compactness.

\paragraph{Caveats.}

The argument depends on two other unverified OpenAI manuscripts: ``Fontaine--Mazur modularity at the prime 2'' and ``A pointwise 2-converse for elliptic curves with rational two-torsion''. The second is needed as written to eliminate Sha[2]. The proof of termination uses a non-constructive compactness argument.

\subsubsection{\texorpdfstring{017 \textperiodcentered{} Irrationality exponent of $\pi$ is 2}{017 · Irrationality exponent of π is 2}}\label{fam:017}
\badges{Lean, 867 files, not in catalogue \badgesep{} combination bordering new \badgesep{} logarithm extension requires a separate argument}

For every $\nu>2$ and all large q, $|\pi-p/q|\ge q^{-\nu}$; hence the Flint--Hills series converges. The previous bound cited by the initial extraction is 7.103 (Zeilberger--Zudilin 2020). The threshold is ineffective.

\paragraph{Extracted instruments.}

\begin{itemize}
\item \textbf{Separated-weight interpolation} (Thm 2.1). Fix integers $m,K\ge1$ and positive rational $w_0,v_0,\theta$ with $0<\theta<1$ and
\[
 K\theta^m<1,\qquad K\frac{w_0}{v_0}\theta^m<1.
\]
Choose positive rational coordinate weights $w_1,\ldots,w_m$ successively above the source's thresholds; each threshold depends only on the fixed data and previously chosen weights. Put $W=(w_0,w_1,\ldots,w_m)$ and $V=(v_0,w_1/\theta,\ldots,w_m/\theta)$. At centers $a_j=(1,c_{j1},\ldots,c_{jm})$, $0\le j<K$, require $c_{0i},\ldots,c_{K-1,i}$ to be pairwise distinct for each coordinate $i$. Polynomials of $W$-degree at most $H$ then realize all packets of $V$-weight strictly less than $H$ under the jets $Y=1+t$, $X_i=c_{ji}+u_i+\log(1+t)$, for all sufficiently large, sufficiently divisible integers $H$. Divisibility depends only on the weights; the lower threshold for $H$ may depend on the centers. The proof uses a Roth/Dyson-type curve inequality, nefness on a weighted blow-up, ampleness, Serre vanishing and Rees pushdown. Related interpolation mechanisms occur in Masser 1981, Fischler 2005 and Fischler--Nakamaye 2014.
\item \textbf{Weighted transverse multiplicity} (Lemma 2.2) and normal-basis rigidity. A B\'ezout-type bound determines the normal-basis index set. The proof combines Philippon's derivative chain with multiplicity estimates.
\item \textbf{Two ways to make the determinant small} (\S{}3). Rows evaluated at rational approximations to periods, $2ijp_i/q_i$, are expressed exactly as test rows of the entire functions $[u^a]P(e^z,z+u)$. The conversion costs $q_i^{-\nu}$ for each unit of transverse index. Rows with the same index test the same functions, allowing Laurent\textquoteright{}s Taylor-degree collision argument to produce a saving that increases with the dimension m. The parameter comparison is feasible exactly when $\nu>2$.
\end{itemize}

\paragraph{Reused.}

Laurent's interpolation determinants; Roth's separated denominators; Philippon zero estimates; Demailly's jets; Mondal's weighted Bézout; Stacks Project.

\paragraph{Caveats.}

An analogue with centers $(r^j,jp_i/q_i)$ for $r\in\mathbb Q_{>0}\setminus\{1\}$ is a proposed extension. The claim $\mu(\log r)=2$ needs a separate proof, including the interpolation hypotheses, rational-height costs and determinant estimates. It is not established by the source comparison in this survey. The initial extraction also flags the interpolation argument's nearest prior art as requiring closer comparison.

\subsubsection{\texorpdfstring{032 \textperiodcentered{} Rational Hodge conjecture for CM abelian varieties}{032 · Rational Hodge conjecture for CM abelian varieties}}\label{fam:032}
\badges{not formalized \badgesep{} combination}

For every complex CM abelian variety the cycle map onto rational Hodge classes is surjective in every codimension; with Milne's results this gives the Tate conjecture for abelian varieties over finite fields and the Hodge standard conjecture for abelian varieties in all characteristics. The repository describes this as one of two results produced outside its fixed computation procedure. The withdrawn eightfolds paper is not cited. The other papers in this family use this result as an input, without supplying inputs to its proof.

\paragraph{Extracted instruments.}

\begin{itemize}
\item \textbf{Four-factor switch} with algebraic contraction. This reduces every CM Hodge class to one tensor relation $v_1+v_2=v_3+v_4$ among odd sign functions, through chains of two-row swaps composed with ample powers and scalar descent. The paper presents this as incremental over Hazama 2003, André--Milne and Gao--Ullmo.
\item \textbf{Surface criterion with fixed-input variation of a rank-two theta tensor.} On a compact ball quotient of SU(V), V a Hermitian 3-space of signature (2,1), Gaussian-completed Taylor coefficients of theta series give holomorphic one-forms and a (2,0)-form $\Xi_W=\sum\Theta_1\wedge\Theta_2$. Freeze the finite theta input and the level, move only the positive-Lagrangian data at the compact place; the period $\int\xi\wedge\bar\xi'$ stays nonzero by continuity, and weak approximation lands on a rational endpoint with the other sign pattern. The initial extraction found no precedent for this step; the rest extends Kudla--Millson and Bergeron--Millson--Moeglin.
\item \textbf{Degree-one CM extraction}: Hecke operators at a split prime satisfy $\Phi^3-D_1\Phi^2+qD_2\Phi-q^3D_3=0$ on the Albanese of the full PEL moduli space; the q-adic valuation of the active eigenvalue reads off the CM sign, and weak admissibility converts slope into Hodge type, so one gets rational maps to the prescribed CM abelian varieties without Rogawski's endoscopic classification.
\end{itemize}

\paragraph{Reused.}

Lefschetz (1,1), Galois trace descent, Hodge--Riemann, Borel--Harish-Chandra, Mumford GIT, Néron models, crystalline comparison, Milne's Hodge-to-Tate theorems, Hazama and Abdulali for the generalized conjecture.

\paragraph{Caveats.}

The source contains about 3,900 lines of unformalized TeX. The Satake normalization, covariance of the constant line and slope formula require consistent sign conventions. The continuity assertion is a central analytic step to check.

\subsection{Analysis and geometry}

\subsubsection{\texorpdfstring{074 \textperiodcentered{} Kakeya maximal conjecture in $\mathbb{R}^{3}$; Kakeya sets in $\mathbb{R}^{4}$}{074 · Kakeya maximal conjecture in ℝ³; Kakeya sets in ℝ⁴}}\label{fam:074}
\badges{not formalized \badgesep{} combination leaning new (3D) \badgesep{} new framework (4D)}

The sources claim $\|K_\delta f\|_{L^3(S^2)}\le C_\varepsilon\delta^{-\varepsilon}\|f\|_{L^3}$ in $\mathbb{R}^{3}$ (the maximal conjecture, stronger than the Wang--Zahl set conjecture), and every set in $\mathbb{R}^{4}$ containing a unit segment in every direction has Hausdorff dimension 4, without measurability or compactness assumptions. The manuscripts have 97 and 175 pages. The earlier bound cited by the initial extraction in $\mathbb{R}^{4}$ is 3.059.

\paragraph{Extracted instruments.}

\begin{itemize}
\item \textbf{Deformed critical-exponent surface with BV selection.} A two-parameter family of critical inequalities balancing a smoothed temporal-deficit cost against a plank-clustering exponent defines $h(p,z)$, the least admissible exponent; the theorem is $h(p,0)=0$. If it fails, $h-z$ is coordinatewise monotone, hence locally BV. BV theory gives approximate differentiability almost everywhere. The source then uses coordinatewise monotonicity to sandwich values at an arbitrary nearby point between comparison values in the density-one set, obtaining ordinary total differentiability almost everywhere; at a differentiability point the partials prescribe the deficit and plank shape of ``equality sequences'', and a minimum-delay argument forces a canonical profile and self-similar ``stationary packets''. This replaces the sticky/plany/grainy induction on scales of Wang--Zahl. The initial search found no matching construction in the Kakeya literature; the familiar ``extremal configuration'' heuristic is the ancestor.
\item \textbf{Two-frame entropy transfer of projection theorems.} Four coordinates of a line refine at rates $(0,\tau,\ell,1)$; two marked events give two frames related by an explicit linear map; a robust planar projection lemma (Ren--Wang Furstenberg in shaded form) plus product-slope pinning shows the product of marginals hits the second frame's support with probability $\le\delta^\kappa$ while the true coupling hits with probability 1; data processing gives $I\ge\kappa\zeta\log D$, contradicting a chain-rule bound. The initial comparison found no earlier use of this information bound to apply a projection theorem to dependent events.
\item \textbf{Chart/horizon framework with narrowness trichotomy} (4D). Lines carry one hidden coordinate; a ``chart'' is a box-times-interval where one quadratic test fits well; the ``horizon'' is the infimum time cost of true terminal systems; charts normalize and compose. Three regimes (potential/conservative field with ``degree restoration'', a Heisenberg-type ``horizontal model'', whole-interval localization) each manufacture a cheaper fit, contradicting extremality. This replaces polynomial Wolff axioms: the quality of quadratic fits becomes the variational quantity.
\end{itemize}

\paragraph{Reused.}

GWZ 2026 convex-clustering estimate (reformulated as a weighted plank lemma, the real dependence of 4D on 3D); Ren--Wang, Wang--Wu; multilinear Kakeya; Hochman--Shmerkin entropy across scales; Orponen--Shmerkin; Christ's quasiextremals; Basu--Pollack--Roy.

\paragraph{Caveats.}

The two manuscripts total 272 pages and use nested limits, slow diagonal choices, triangular subsequences and almost-everywhere offsets. Their order of limits needs detailed checking beyond the selected comparisons made here. The three-dimensional argument cites another unformalized OpenAI restriction paper as an alternative route. The initial extraction identifies Rios--Sawyer as an uncited antecedent. Literature comparisons remain incomplete.

\subsubsection{\texorpdfstring{087 \textperiodcentered{} The Mahler conjectures}{087 · The Mahler conjectures}}\label{fam:087}
\badges{original inspection: symmetric Lean coverage only; broader solution modules supplied, unbuilt \badgesep{} new identity}

For origin-symmetric convex bodies $K\subset\mathbb R^n$, the source claims $|K||K^\circ|\ge4^n/n!$, with equality for linear images of Hanner polytopes. For general convex bodies, use the Santaló point $s(K)\in\operatorname{int}K$ and define
\[
 P(K)=|K|\,|(K-s(K))^\circ|.
\]
The general claim is $P(K)\ge(n+1)^{n+1}/(n!)^2$, with equality exactly for simplices; for a simplex its Santaló point is its centroid. Thus the equality statement for the polar is made after the indicated centering. In the symmetric case $s(K)=0$. The third source claims $c_G(\operatorname{int}K\times\operatorname{int}K^\circ)=4$ for origin-symmetric $K\subset\mathbb R^n$ and $n\ge2$, with embeddings at every capacity strictly below 4. The initial extraction inspected symmetric formalization coverage only. The snapshot also supplies general-Mahler and symmetric-polar solution modules; a \texttt{sorry} in a comparator challenge is a statement placeholder, not evidence that the separate solution module is a stub. None was built here.

\paragraph{Extracted instruments.}

\begin{itemize}
\item \textbf{The conformal lens.} $F(z)=\frac8{\pi^2}\sum(-1)^jz^{2j+1}/(2j+1)^2$ maps the disk onto a convex lens whose boundary, under uniform angle, has uniform imaginary coordinate and a fair-sign real coordinate (Lemma 3.1, ``lens boundary law''). For a strip body $A=\{|b_i\cdot X|\le1\}$ form $\Omega=\{z:\ b_jz\in D\}$. The map is Gross's 2019 conformal Skorokhod example, as the paper acknowledges. The initial comparison found no earlier application to volume products. The lens replaces the tube domain used by Nazarov and Berndtsson.
\item \textbf{Holomorphic mass pushed to the boundary.} For a holomorphic map $f:\mathcal U\subset\mathbb C^d\to\mathbb C^N$, $N\ge d$, assume that $f^{-1}(0)=\{0\}$, its leading homogeneous part has degree $k\ge1$ and is nonzero away from zero, and every sublevel $\{|f|^2\le s\}$, $0<s<1$, is compact in $\mathcal U$. The integral of the sum of squared $d$-by-$d$ Jacobian minors over $\{|f|^2<1\}$ is then at least $(\pi k)^d/d!$ (the source's special case of Lelong comparison). Apply to $f_k=(g(b_1z)^k,\dots)$ with $g=F^{-1}$, substitute $u_i=\rho_i^{1/k}e^{i\theta_i}$, let $k\to\infty$: Fatou pushes the mass onto boundary samples, giving $S=\sum_IP_I\ge1$ where $P_I$ is the probability that d independent uniform boundary samples solve all strip constraints. The rescaling that turns Lelong mass into feasibility probabilities had no precedent in the literature examined.
\item \textbf{Feasible-simplex identity} (Prop. 5.1): $1\le S=\frac{d!}{4^d}\int_A|\Sigma_X|\,dX\le\frac{d!}{4^d}|A||A^\circ|$, where $\Sigma_X$ is a union of interior-disjoint simplices $\mathrm{conv}(0,\varepsilon_ib_i)$ over feasible signed bases, with an explicit deficit. Equality forces ``lens-cost medians'' for every triple, then Lima's criterion and Hansen--Lima give the Hanner classification. The identity is classified as new in the initial extraction.
\end{itemize}

\paragraph{Reused.}

Gross 2019; Demailly; Cauchy--Binet; Hanner's formulas; Lima and Hansen--Lima (the Lean reproves the median-to-Hanner step itself); Fradelizi--Meyer; Klartag 2018 for the general case.

\paragraph{Caveats.}

The general-body argument uses about twenty selected decimal constants and SymPy certificates. Its reasoning summary leaves the conclusion dependent on uniform scalar inequalities. The symmetric paper presents a self-contained argument without that qualification.

\subsubsection{\texorpdfstring{362 \textperiodcentered{} Relativistic Vlasov--Maxwell, 3D, large data}{362 · Relativistic Vlasov–Maxwell, 3D, large data}}\label{fam:362}
\badges{Lean, 207 files, continuation step included \badgesep{} combination}

For the one-species system with no background charge, the source assumes
\[
\begin{gathered}
 f_0\in C_c^\infty(\mathbb R_x^3\times\mathbb R_v^3),\qquad f_0\ge0,\\
 E_0,B_0\in C_b^\infty(\mathbb R^3;\mathbb R^3)\cap L^2(\mathbb R^3;\mathbb R^3),\\
 \nabla\cdot E_0=\rho_{f_0},\qquad \nabla\cdot B_0=0,\qquad
 \rho_{f_0}(x)=\int_{\mathbb R^3}f_0(x,v)\,dv.
\end{gathered}
\]
Here every spatial derivative of the fields is bounded, the fields have finite energy, and the nonnegative phase-space datum is smooth and compactly supported. Under these hypotheses the paper claims a unique global classical solution, smooth on each finite time interval and with compact phase-space support there, without a size restriction. The repository supplies a formalization covering a continuation statement; it was not compiled or kernel-checked for this survey.

\paragraph{Extracted instruments.}

\begin{itemize}
\item \textbf{Signed momentum increment bootstrap} (Prop. 2.1): for a receiver with energy $w\le P$, $|V(t_2)-V(t_1)|\le MP\sqrt{t_2-t_1}+A(t_2-t_1)P^2L/w$, the norm ``taken after integration''. With $P=2^n$ each doubling takes time $\ge c/\log 2^n$; $\sum1/n$ diverges, so the momentum support stays bounded on every finite interval and the Glassey--Strauss continuation theorem completes the argument. This replaces pointwise field bounds that lose $\sqrt w$. The time-integrated increment idea is in Cheverry--Ibrahim 2025; the initial extraction identifies the bootstrap formulation and logarithmic doubling argument as the new features.
\item \textbf{Signed force identity along source labels} (Lemma 6.1): the Liénard--Wiechert kernel along a source's retarded branch equals a total derivative plus residuals with no source acceleration and no inverse powers of the source deficit, at the cost of a term proportional to the receiver's own acceleration. The physical ancestor is the Heaviside--Feynman form; the analytic version with moving receiver is new in this setting.
\item \textbf{Occupation of retarded bins, direction-change counting, selection} (\S{}\S{}4, 7, 8): at most $CL^6/\delta$ direction changes of size $\delta$; an improved occupation bound in a window; a selection weight with $\varepsilon=10^{-5}$ and an exponent contradiction (0.496z versus 0.41z) that makes the receiver coefficient $\le Cw^{-1/2}$, cancelling the loss exactly. The occupation estimate adapts Pfaffelmoser's method to the light cone; the counting and closure are classified as new.
\end{itemize}

\paragraph{Reused.}

Glassey--Strauss representation; Luk--Strain cone flux; Lazarovici label change; Pallard's cone Jacobians; Pfaffelmoser; Kirchhoff; Gagliardo--Nirenberg; Picard.

\paragraph{Caveats.}

The reasoning summary states that the cited continuation theorem alone does not establish the required estimates, whereas the paper presents them as proved. Cheverry--Ibrahim is cited six times in the summary and omitted from the paper. The exponent comparison has a small margin. A repository formalization is supplied, but no machine check was run for this survey. Uniqueness is asserted only among global classical solutions.

\subsubsection{\texorpdfstring{304 \textperiodcentered{} Hilbert--Smith conjecture in every dimension}{304 · Hilbert–Smith conjecture in every dimension}}\label{fam:304}
\badges{not formalized \badgesep{} new}

Every locally compact second-countable Hausdorff group acting faithfully and jointly continuously on a connected Hausdorff second-countable topological manifold of positive finite dimension without boundary is a Lie group. The equivalent exclusion concerns faithful jointly continuous $\mathbb{Z}_p$ actions on such manifolds, for every prime $p$. Pardon's 2013 proof covered dimension 3 via Smith theory and geometrization; the surveyed argument does not use Yang's theorem, Smith theory, cohomological dimension or cell-like maps.

\paragraph{Extracted instruments.}

\begin{itemize}
\item \textbf{Compact-source test categories} $T(X)$: the thick subcategory of $D^b(\mathrm{Sh}_R(X))$ generated by $Rf_*R_B$ for arbitrary continuous maps f from compact polyhedra. Stalks may be infinite-dimensional and images non-constructible, yet Verdier duality restricts with invertible biduality, $Rg_!=Rg_*$ is duality-preserving for any continuous g, and Verdier quotients have a ``germ description'' via a cutoff lemma. This replaces Woolf's perfect/constructible Witt groups (whose excision may fail) and Ranicki--Weiss's chain-level duality (the exact antecedent of the cutoff).
\item \textbf{Fixed sphere lattice / bounded comparison.} The semialgebraic-to-continuous comparison $F_d(S^d)\to E_d(S^d)$ has kernel and cokernel killed by $2^{e_d}$, via Balmer--Walter localization, a dense-doubling map with $\iota d=d\iota=2$ (classical exponent-2 bound of Hornbostel--Schlichting), and Hardt triviality supplying an integer-valued signature. So $E_d(S^d)\otimes\mathbb Q=\mathbb Q\bar u$ with image in $L_d^{-1}\mathbb Z\bar u$ for $L_d$ fixed before any action data: ``integral before rationalization''.
\item \textbf{Character classes and equal parts in a fixed lattice.} The orbit sheaf splits into eigensheaves of the $\mathbb{Z}_p$ generator; finite partitions give self-adjoint projectors becoming germ endomorphisms; for $m=p^k$ faithfulness yields a region over which the orbit map is an m-sheeted trivial cover; a locally constant unit multiplier cyclically permutes the parts as a real isometry; Serre's odd-sphere degree maps concentrate the class; each part equals $(4/p^k)\bar u$, impossible once $p^k\gt4L_d$. The initial search found no earlier version.
\end{itemize}

\paragraph{Reused.}

Newman 1931 and Pardon's chart reduction; Kashiwara--Schapira duality; Balmer's triangular Witt groups; Woolf's cone lemma; Ranicki--Weiss cutoff; Hardt triviality; Dranishnikov--Ferry--Weinberger telescope; \v{C}ech theory and Haar averaging.

\paragraph{Caveats.}

Several steps require closer examination: whether the Balmer--Walter hypotheses pass to localizations of the generated categories; whether $Rt_*$ commutes with the infinite eigensheaf sums used; and whether the restricted pairing is isometric in the germ quotient. The odd stabilization $d>n$ is essential to the argument.

\subsection{Combinatorics and theoretical computer science}

\subsubsection{\texorpdfstring{102 \textperiodcentered{} Unique Games Conjecture}{102 · Unique Games Conjecture}}\label{fam:102}
\badges{Lean, 475 files, full statement, PCP theorem reproved \badgesep{} new gadget}

For every fixed $\varepsilon,\delta\in(0,1/2)$ there is a fixed alphabet $\F_2^s$ and a deterministic polynomial-time reduction from 3SAT to bipartite Unique Games with translation constraints, value $\ge 1-\varepsilon$ on satisfiable and $\le \delta$ on unsatisfiable inputs. The formal object fixes alphabet, machine and predicates; $\varepsilon$ and $\delta$ are independent. The development also formalizes a Dinur-style PCP theorem, Håstad's 3LIN and the Khot--Minzer--Safra Grassmann expansion theorem from scratch.

\paragraph{Extracted instruments.}

\begin{itemize}
\item \textbf{Latent alphabet gadget} (Lemma 2.1, \S{}2). The construction consists of a binary space $V\supseteq K=\F_2^\ell$, a nonlinear map $C:V\to K$ with $C(x+k)=C(x)+k$, and a rational noise law $\mu$ on V with $\Pr[C(x+a)\ne C(x)]\le p^*$ for independent uniform $x\in V$ and $a\sim\mu$, while every subspace of $V^*$ with large restricted rank detects $a\sim\mu$ with probability $\ge 1/8$, the rank threshold depending only on $p^*$, fixed before $\ell$. It is built from the quadratic block $Q(x)=(x_2x_3,x_1x_3,x_1x_2)$ over $\F_{2^d}$, update spaces indexed by the $q^2+q+1$ lines of $\F^3$ (the quadratic term cancels on a line), a trace-duality ``alignment equation'' $z\cdot g(z)=\sqrt{z_1z_2z_3}$ characterizing rank loss, recursive concatenation with a harmonic potential bounding rank loss by 4$\theta$ per level. The resulting noise rarely changes the honest evaluation (completeness $\approx 1$) yet forces 0.9 of Fourier mass below bounded rank, so the existing degree-two KMS inverse theorem applies unchanged. Barak--Kothari--Steurer 2018 proposed higher-degree tests for exactly this and could not compose them; here degree stays two and the noise changes. The initial search found no matching gadget.
\item \textbf{Exact table keys with translation-orbit folding} (\S{}3.3): tables are affine maps on homogeneous answer coordinates, keys invariant under sparse projections, vertices are translation orbits, so the sparsely projected prover shares the UG labels ``without an additional consistency test''. Combination of Khot--Saket virtual tables and Minzer--Zheng.
\item \textbf{Row-fiber abundance by erasure} (Lemma 4.4): KMS gives one good slice per fiber; re-randomizing the function where good fibers are sparse and a Hoeffding union bound over slice descriptions shows good advice has probability $\ge\eta2^{-\ell r_s}$ independent of k. The initial extraction identifies this as the new step in the advice game; the reasoning summary records ``KMS just existence of one fiber, need union randomization argument''.
\end{itemize}

\paragraph{Reused.}

Håstad's 3LIN gap; KMS Grassmann expansion through BKS's matrix chart; Khot--Saket folding and list decoding; KMS/DKKMS/Minzer--Zheng smooth advice game; Dinur--Steurer parallel repetition twice; Schwartz--Zippel. Not used: Majority is Stablest, high-dimensional expanders, sum-of-squares.

\paragraph{Caveats.}

The argument does not use other OpenAI manuscripts as inputs. Its reduction exponent and alphabet size are tower-like in $1/p^*$, which accounts for compatibility with the Arora--Barak--Steurer subexponential algorithm, although the paper does not discuss this point. The reasoning summary records that UGC arose as a subgoal of a different prompt about the basic SDP threshold.

\subsubsection{\texorpdfstring{107 \textperiodcentered{} Matrix multiplication exponent $\omega \le 9/4$}{107 · Matrix multiplication exponent ω ≤ 9/4}}\label{fam:107}
\badges{Lean, straight-line definition of $\omega$ over $\mathbb{C}$ \badgesep{} new as a proof of an $\omega$ bound}

$\omega(\mathbb C)\le9/4$, compared with the previous bound 2.371177 cited by the initial extraction. The 13-page argument is non-constructive: it bounds every point of Strassen's asymptotic spectrum and gives ``no competitive finite matrix size''.

\paragraph{Extracted instruments.}

\begin{itemize}
\item \textbf{Finite separation of a shared leg and the shared-leg entropy inequality} (Prop. 3.1, Cor. 3.2). Starting with $A=\sum_hA_h$ and a shared X-leg, the construction takes $L=5M$ copies and substitutes roots of unity. The resulting r-sum is a DFT that retains the terms satisfying $u-v+2(g-h)=0$. On these terms, the degeneration weights $g^2,\ hu-h^2,\ -hv$ add to $(g-h)^2$, so the weight-0 component is the complete direct sum $\bigoplus_hA_h\otimes B_X(M)$. For nonzero tensors $T_i$ sharing the first space and occupying matched, pairwise disjoint sectors on the other two legs, every spectral point $\lambda$ and probability vector $q$ satisfy $\lambda(T)\ge e^{p_XH(q)}\prod\lambda(T_i)^{q_i}$, where $T=\sum_iT_i$ and $H(q)=-\sum_iq_i\log q_i$. The construction retains all blocks at a factor-5 cost; Salem--Spencer hashing and laser-method zeroing-out retain at most one block per shared variable. Its ingredients are the Cohn--Umans DFT, Strassen degeneration, and the method of types. The initial literature comparison identified no matching combination.
\item \textbf{Spectral growth on polynomial multiplication.} $C(a,b)=\sum x_iy_jz_{i+j}$ has rank $a+b-1$; symmetrize over $S_3$ to a profile $P(a,b)\le(a+b-1)^{1/t}$. Two degenerations fed into the inequality (a Clebsch--Gordan filtration giving discrete concavity in b, and a three-window partition giving tripling) force $P(a,a)\ge a^{4/3}$, hence $3t\le9/4$. Strassen duality converts the uniform spectral bound into an asymptotic-rank bound. The source compares its spectral argument with barriers for specified tensor-extraction methods, including the laser and universal methods. Its proof constructs no degeneration into matrix products; this comparison does not establish that it evades every possible barrier. Alman--Li's 2026 ``Strassen calculus'' is the uncited template, used there only for $\mathrm{cw}_q$-type bounds.
\item \textbf{Detecting characters} (App. A): a fixed-point (Schauder--Tychonoff) proof of the needed case of Strassen's spectral theorem. The initial extraction classifies this reproof as incremental.
\end{itemize}

\paragraph{Reused.}

Bini interpolation; Strassen's asymptotic spectrum; combinatorial degeneration; method of types; Clebsch--Gordan exact sequence.

\paragraph{Caveats.}

The stated theorem concerns $\mathbb C$. The argument bounds spectral profiles through a specified collection of inequalities; a better exponent would require a further argument. No Lean build was performed for this survey.

\subsubsection{\texorpdfstring{157 \textperiodcentered{} Counterexample to Hadwiger's conjecture}{157 · Counterexample to Hadwiger's conjecture}}\label{fam:157}
\badges{Lean, 276 files, no explicit graph \badgesep{} new}

The source constructs arbitrarily large m-vertex graphs with independence number $\le 2$ and connected-matching number $<m/100$, where every pair of matching edges must have an edge between their endpoint sets. The resulting Hadwiger number is $h(G)<26m/75+2/3$, hence $h(G)<0.347m$ for sufficiently large $m$, while $\chi\ge\chi_f\ge m/2$. This would also refute the Reed--Seymour fractional weakening and the Füredi--Gyárfás--Simonyi connected-matching conjecture. Parameters: $g=10^9+1$, $m=2^{1000\cdot g\cdot M_0\cdot n}$; the Lean extracts m from an eventually-true filter statement.

\paragraph{Extracted instruments.}

\begin{itemize}
\item \textbf{The hole relation with tensor certification} (\S{}2). Vertices are random injective frame maps over $\F_2$ with prescribed self-Gram; two vertices have a ``hole'' iff witnesses $\lambda_i$, $\lambda_j$ in a cut-profile space satisfy a shared-tensor equation and gradient identities $u_jU_i=r(\lambda_i)$, $u_iU_j=r(\lambda_j)$ with $r(x)=a+T(x,\cdot)$, T symmetric bilinear, and $a(\lambda_i)+a(\lambda_j)=1$. Summing the gradient identities around a hypothetical triangle cancels all bilinear terms and leaves $0=1$, so holes are triangle-free deterministically and $\alpha \le 2$ follows directly. This replaces the pseudorandom triangle-free graphs of prior $\alpha=2$ attacks. The initial comparison found no matching construction.
\item \textbf{Entropy-minimizing container argument} (\S{}3): Kleitman--Winston fingerprinting on the conflict graph of units whose progress measure is the minimum relative entropy of a law obeying marginal caps; each recording raises the minimum by $-\log(1-\varepsilon)$; when no feasible law remains, a Ford--Fulkerson cut certifies the residual is covered by small sets. It transfers a distributional supersaturation theorem (``any capped unit law has four-cross-hole probability $\ge2^{-100gN}$'') to the random sample. The initial search found no precedent for this potential or the terminal min-cut certificate.
\item \textbf{Scalar recipe / gradient realization / collision machine} (\S{}\S{}4--13): witnesses as sums of $\le 7$ atoms, keys forcing all four holes from O(n) scalar Gram bits, overlap forced by min-entropy peeling, Walsh bounds, a Frieze--Kannan ``twisted cut norm'' increment, and a characteristic-2 Boolean flat-extension lemma. This is an extensive collection of problem-specific constructions; the initial extraction did not settle its novelty.
\end{itemize}

\paragraph{Reused.}

Plummer--Stiebitz--Toft $\alpha \le 2$ reduction; Cambie's connected-matching-to-minor count; Kleitman--Winston; Csiszár; Ford--Fulkerson; weak regularity; Laurent--Mourrain flat extension.

\paragraph{Caveats.}

The argument gives neither an explicit graph nor an explicit clique-minor number. Many subsequence/compactness steps are non-quantitative by design (``No exponential rate uniform over all stages is asserted''). Illingworth--Steiner's separate result of 28 September 2026 disproves the dominating variant~\cite{src:dominating}, which requires a stronger kind of clique-minor model. A counterexample to that variant need not violate Hadwiger's conjecture itself. Both constructions use graphs with $\alpha\le2$. The repository supplies a Lean development for the claim surveyed here, but it was not compiled for this survey.

\subsubsection{\texorpdfstring{158 \textperiodcentered{} The plane is not 5-colorable}{158 · The plane is not 5-colorable}}\label{fam:158}
\badges{Lean, unrestricted statement, $\approx 24{,}500$ lines \badgesep{} new}

Every 5-coloring of $\mathbb{R}^{2}$ with arbitrary color classes has a monochromatic unit pair, so the chromatic number of the plane is 6 or 7, in ZFC. The formal statement is literally $\neg\exists\, c:\mathbb{C}\to\mathrm{Fin}\,5$ proper. The proof uses neither a SAT solver nor a vertex list; its argument is structural and ergodic.

\paragraph{Extracted instruments.}

\begin{itemize}
\item \textbf{Transfer of colorability via Haar rigidity of wild character laws} (\S{}\S{}2--4). Let E be the algebraic plane, K its algebraic unit circle, D the dual of E with discrete topology, $C\subset D$ the continuous characters. A K-invariant probability on D is ``wild'' if it gives C measure zero. The stated theorem is that every wild K-invariant probability is Haar. Its proof combines a Furstenberg--Zimmer compact tower with an explicit approximation property, a transfinite relative-singularity lemma, conditional multiple averages with van der Corput, a radial difference inequality $w(|A-B|)\le w(A)+w(B)$, and a number-field general-position lemma. To obtain the coloring consequence, Følner-average any proper coloring over $E\rtimes K$; the spectral ``wild'' part of each color indicator is uncorrelated at all algebraic translations and can be discarded, leaving a measurable ``weak'' coloring with almost no monochromatic unit pairs. Thus a proper k-coloring exists if and only if a weak measurable k-coloring exists, for every k. The literature since Falconer treats measurable and unrestricted chromatic numbers as possibly different; the claimed transfer equates existence in these two formulations. The initial literature comparison found no matching transfer theorem.
\item \textbf{Angular palettes, binary centers, transition graph} (\S{}\S{}5--6): weak-$*$ limits of circle samples around a center; a 60$^\circ$ exclusion; the bound ``at most two colors per direction a.e.'' proved by an expanding-action bounded-distortion argument on six circle maps. Combination of Falconer/Townsend/Sokolov--Voronov density mechanics with a new dynamical separation.
\item \textbf{Common exclusion continua} (\S{}7): threshold and renormalize indicators into a continuous map to the 1-skeleton of the 4-simplex; a Brouwer cover obstruction forces a non-null loop; midpoint preimages become continua that each exclude two colors in a crescent. Combination.
\end{itemize}

\paragraph{Reused.}

de Bruijn--Erd\H{o}s; Falconer, Payne; Furstenberg--Zimmer; Bochner/SNAG; Følner means; Fourier decay of the circle; Brouwer; Hatcher covering theory; Moser spindle with a rational placement certificate.

\paragraph{Caveats.}

The rigidity theorem is dense and relies on a transfinite induction; the initial extraction places it inside the Lean closure (Rigidity.lean). The inspected comparator configuration disables the optional nanoda independent kernel recheck; no build was run here. The paper never exhibits a finite 6-chromatic graph, consistent with the method. A footnote disputes a pseudonymous manuscript; that dispute is not assessed here.

\subsubsection{\texorpdfstring{159 \textperiodcentered{} Erd\H{o}s's reciprocal-sum conjecture}{159 · Erd\H{o}s's reciprocal-sum conjecture}}\label{fam:159}
\badges{original label: companion only; repository also links a headline development, unbuilt \badgesep{} new architecture}

The source claims $r_k(N)\le C_kN\exp(-c_k(\log N)^{\varepsilon_k})$ for each fixed $k\ge3$, with $c_k,\varepsilon_k>0$, and hence progressions of every finite length in any set with divergent reciprocal sum. The dyadic deduction is valid: for a $k$-progression-free set $A$,
\[
 \sum_{n\in A\cap[2^j,2^{j+1})}\frac1n
 \le 2^{-j}r_k(2^{j+1})
 \ll_k \exp(-c'_k j^{\varepsilon_k}),
\]
and the last expression is summable in $j$. By contrast, $\sum_{m\ge2}\exp(-(\log m)^c)$ diverges for $0<c\le1$ and converges for $c>1$. The original discussion of the weaker Leng--Sah--Sawhney log-log bound must be understood in the former range; divergence is not true for every positive $c$. The 198-page source uses the LSS inverse theorem as an external input.

\paragraph{Extracted instruments.}

\begin{itemize}
\item \textbf{Triangular polynomial cells with descending precision budgets.} A density-increment region is a nested system $\|b_h-C_h(u,b_{\lt h})-l_h\|_\infty\le w_h$, $h=1,\dots,D$, where integer blocks of weight h are determined exactly, $C_h$ is a weighted-degree polynomial with unbounded coefficients, and $Q_h=\log(2/w_h)$ is the precision. The key inequality is that the loss at weight i, $Q_i'-Q_i$, is bounded by a polynomial in the precisions at weights above i, never in $Q_i$ itself. With D fixed, descending induction keeps all precisions polynomial in $p=\log(1/\alpha)$ over O(p) increments, which is exactly a quasipolynomial bound. This replaces increments that reset to a sub-progression and lose a power each round. The trace shows it emerging in Part II as the fix for a compounding-resources recurrence. The initial search found no earlier use of this terminology or architecture.
\item \textbf{Shift comparison and positive counting}: a one-sided structure theorem at higher nilsequence degree (degree 1 is Kelley--Meka sifting plus Schoen--Sisask), with a two-kernel trick converting a deficit into a surplus detectable by upper bounds. The initial comparison identified binary-system and corners analogues, and classified the stated higher-degree propagation as new.
\item \textbf{Relative lifting with additive rank and exact marked projection}: sample integer-affine paths inside a cell, recover an ambient nilsequence from local ones, keep every old integer slot exact; returned rank independent of ambient dimension. This is a quantitative local-to-global transfer, classified as a new formulation; its correctness was not established for this survey.
\end{itemize}

\paragraph{Reused.}

LSS inverse theorem; Schoen--Sisask almost-periodicity; Kelley--Meka sifting; Green--Tao polynomial orbits; Leng's equidistribution; Candela--Sisask densification; CRT prime boxes.

\paragraph{Caveats.}

The quantitative upper bound was not verified here. The repository documentation links a development for the reciprocal-sum consequence, beyond the lower-bound companion noted in the initial extraction; that development was not built here. The paper claims an improvement from log log to log for every k at once. The trace's Part I never closed (``circularity concern unresolved'') and Part II ends in ``claimed/asserted''. The argument has a long chain of dependencies. An appendix extends LSS to product groups of dimension $(2+p)^A$ with a step-drop lemma uniform in parameter count, which the model itself flagged as delicate.

\subsection{Algebra and operator algebras}

\subsubsection{\texorpdfstring{196 \textperiodcentered{} Counterexample to Kaplansky's zero-divisor conjecture}{196 · Counterexample to Kaplansky's zero-divisor conjecture}}\label{fam:196}
\badges{Lean, 124 files, includes finite 2-dim K(G,1) \badgesep{} combination; separator step new}

The source constructs a finitely presented torsion-free group with a finite two-dimensional classifying complex and nonzero $\alpha$, $\beta$ in $\F_2[G]$ with $\alpha\beta=0$. The construction is non-constructive: 8258 generators, a random graph on $n\gtrsim e^{300000}$ vertices. The initial extraction cites 2026 revisions by Mineyev and Shin that describe the conjecture as open.

\paragraph{Extracted instruments.}

\begin{itemize}
\item \textbf{Odd cross-intersection parity mechanism.} The alphabet consists of the 16513 points of PG(2,128) plus three letters; each vertex of two graphs gets an outgoing label set whose ordinary part is a projective line, arranged so every pair of label sets across the two graphs meets in an odd number of points. Cone off components; $\alpha$ and $\beta$ are sums of root-component path labels. Stepping both coordinates by a common letter leaves the product element invariant, every vertex of the product graph has odd degree, so every component has even size and every coefficient of $\alpha\beta$ vanishes mod 2, ``independently of how edges will be matched''. This relaxes Mineyev's even-product matching and the Rips--Segev non-unique-product requirement to ``odd-degree product graph'' realized by an oddtown-style design.
\item \textbf{Bounded-pattern estimate for a girth-conditioned random-matching model with prescribed types} (Prop. 3.1): conditional on girth at least $L$, the probability that such a system exists is $o(1)$ as the admissible graph size tends to infinity: it has at most $K$ positive-length paths, total length $L\le H\le CL$, at most $I$ interval pairs, and at most $\varepsilon H$ unpaired occurrences. All paths are cyclically immersed and closed, except possibly one immersed interval starting at a prescribed root. Paired occurrence intervals are disjoint, obey the signed-label matching rules, and traverse distinct underlying undirected edges. The constant $\varepsilon$ is fixed before $K,C,I$, which remain fixed in the limit. The proof proceeds via multiplicity stages, ``nearly aligned blocks'' from a simultaneous Dirichlet choice, and a transfer-matrix bound. This reverses the Gromov--Ollivier--Osajda model (fixed graph, random labels) into fixed label types, random graph, with no small-cancellation hypothesis.
\item \textbf{Planar-separator extraction plus cone-picture band surgery} (Props. 4.4, 5.1): any reduced spherical arrangement is cut by Lipton--Tarjan recursion to bounded clusters that the estimate excludes; the surgery along bands joining two occurrences of a lifted edge proves asphericity and root protection together. The initial extraction classifies the separator step as new in this setting and the surgery as incremental over Gruber and Bereznyuk.
\end{itemize}

\paragraph{Reused.}

Graphical presentations; projective-plane incidence parity; configuration model with switchings; Ollivier's density counting; Lipton--Tarjan; Brown's periodic-resolution argument.

\paragraph{Caveats.}

The result concerns characteristic two and does not provide a concrete presentation. The source explicitly addresses quantifier order and the dependence introduced by conditioning. Its Lean development was not compiled here. A separate direct-finiteness manuscript dated 4 October adapts this machinery and is not formalized.

\subsubsection{\texorpdfstring{197 \textperiodcentered{} Counterexample to Kaplansky's direct-finiteness conjecture}{197 · Counterexample to Kaplansky's direct-finiteness conjecture}}\label{fam:197}
\badges{Lean, 125 files, search certified at $h \le 217$ \badgesep{} new criterion}

The construction gives a finite field of characteristic two, a finitely presented group G with odd torsion, and a, b in K[G] with $ab=1,\ ba\ne1$. Such a group would be non-sofic by Elek--Szabó, and the construction would also contradict Gottschalk's surjunctivity conjecture.

\paragraph{Extracted instruments.}

\begin{itemize}
\item \textbf{The incidence criterion} (Thm 3.1): let $\ell$ be an odd prime and $K$ a finite field of characteristic two containing a primitive $\ell$th root of unity. Take a finite hypergraph with $t$ points and nonempty blocks, every point on exactly $\ell+1$ blocks, and bipartite incidence graph of girth at least 12, together with vectors $x_v\in K^m$ with $\sum_{v\in f}x_vx_v^\T=0$ per block and $\sum_vx_vx_v^\T=-I_m$, with $0<m<t\le m\ell^2$; these data yield the counterexample. The key identity uses the fact that $\ell+1$ is the number of lines through the origin of $\F_\ell^2$, so local copies of $\F_\ell^2$ glue along the blocks, and $\sum_{u\in K_v}u=1+\sum_{f\ni v}(\sum_{u\in L_f}u-1)$ turns the two identities into $UT=I_m$ after averaging. This gives a local-to-global certificate for failure of direct finiteness; the closest relative is Dykema--Heister--Juschenko's universal left-invertible-element groups, which encode cancellation patterns but give no such criterion.
\item \textbf{Polynomial pairing on a punctured affine space} (Prop. 2.3): Newton interpolation polynomials for a set of holes give the vectors; products of degree $\le q-2$ sum to zero on every full line by the power-sum identity, and removing the holes leaves $-I$. The multiplication property also appears in Golowich--Lin and Dinur--Liu--Zhang codes, which the paper acknowledges; it does not invoke a theorem from those constructions.
\item \textbf{Hole trades with random slabs and the local lemma} (\S{}5): replacing the fan of lines through a hole by parallel lines in a plane preserves every other point's degree; Linnik's theorem supplies $\ell\equiv-1\bmod q^{100}$ so bundles of q directions exist; short incidence cycles are enumerated and excluded by the Lovász local lemma. The resulting regular high-girth line system is classified as a new construction.
\item \textbf{Character-twisted HNN packing and the sandwich certificate} (\S{}4): $\tau_v$ conjugates $K_v$ onto a reference copy while twisting by a linear form, so t coordinates share m slots with orthogonal character idempotents; $J_0\Delta C=e(TU-P)$ lives in the smaller algebra K[H], where the homomorphism ``every group element $\mapsto 1$, $z \mapsto \zeta$'' gives $X^\T X-I_t\ne0$ by rank. The planar Euler-characteristic embedding argument uses a form of small cancellation; the Dykema--Juschenko finite-group step is classical and made explicit.
\end{itemize}

\paragraph{Reused.}

Elek--Szabó for non-soficity; Dykema--Juschenko reduction; HNN embedding; Newton interpolation; Lovász local lemma; Linnik via Xylouris.

\paragraph{Caveats.}

The construction was not executed at the prescribed scale (q up to $2^{217}$, dimension $\ge 1200$); the certificate is algebraic. The Lean statement does not include non-soficity. The repository supplies a formalization for the odd-characteristic companion, but not for the torsion-free companion dated 4 October.

\subsubsection{\texorpdfstring{248 \textperiodcentered{} Thompson's group F is not amenable}{248 · Thompson's group F is not amenable}}\label{fam:248}
\badges{Lean, conceptual, F defined from scratch \badgesep{} new use, leaning new}

F admits no positive normalized left-invariant mean, settling Geoghegan's 1979 question in five pages. The formal statement defines F as dyadic piecewise-linear homeomorphisms and the invariant mean on $\ell^\infty(F)$; it gives no finite certificate, and the constants are existential.

\paragraph{Extracted instruments.}

\begin{itemize}
\item \textbf{Recursive Hilbert-ball-valued coloring of dyadic partitions.} For a dyadic partition T of $[0,1]$, define $p(T)=f\big(\tfrac1D\sum_jp(T_{I_j})\big)$ when T respects D fixed internal basic dyadic intervals with strictly positive pairwise gaps, and define it to be 0 otherwise. Here $T_{I_j}$ is the normalized restriction. The recursion is well founded because restriction decreases the number of cells. Applied to the image partitions $gT^{(n)}$, it produces colors $X_I(g)$ for which the parent color equals f applied to the mean of the descendant colors. Moore\textquoteright{}s 2013 reformulation uses finite, nonrecursive tree colorings. The initial comparison found no earlier version of this recursive rule with a fixed-point-free f.
\item \textbf{Benyamini--Sternfeld displacement map as an amenability obstruction.} Let $f$ be an $L$-Lipschitz self-map of a Hilbert unit ball with $\|f(x)-x\|\ge\delta>0$ everywhere (Benyamini--Sternfeld 1983; the source constructs $\delta=1/2$ on $L^2([0,1];\mathbb R^2)$). For each $g$ in the nonempty finite averaging set $A\subset F$, the recursive identities give the pointwise estimate
\[
 \delta^2\le\|m(g)-f(m(g))\|^2
 \le\frac{L^2}{D}\sum_{i=1}^D\|z_i(g)-m(g)\|^2.
\]
The correlation comparisons control the finite averages, not each summand pointwise. With $\mathbb E_A$ denoting uniform averaging over $g\in A$ and $\eta$ the maximum normalized translation defect for the selected transports,
\[
\begin{aligned}
 \delta^2&\le\frac{L^2}{D}\sum_{i=1}^D\mathbb E_A\|z_i-m\|^2\\
 &\le L^2\left(\frac4D+4\left(1-\frac1D\right)\eta\right).
\end{aligned}
\]
Choosing $D>4L^2/\delta^2$ gives a positive lower bound on $\eta$, contradicting the Følner criterion. The map is non-affine and $F$ does not act on it; the initial extraction found no earlier use of this displacement theorem for this amenability problem.
\item \textbf{Correlation comparison via pair transport}: for internal basic dyadic intervals $I<J$ and $I'<J'$ in $(0,1)$, where $I<J$ means $\max I<\min J$, an element of $F$ carries $I$ affinely onto $I'$ and $J$ affinely onto $J'$. The positive gaps and internal endpoints are part of the hypothesis. This compares averaged color correlations for the strictly separated pairs in the construction; it makes no such assertion for nested pairs or arbitrary boundary-touching intervals.
\end{itemize}

\paragraph{Reused.}

Følner's criterion (formalized forward direction); Cannon--Floyd--Parry partition model; Benyamini--Sternfeld; convexity of the square norm.

\paragraph{Caveats.}

The five-page argument addresses Geoghegan's amenability question. The pointwise and averaged estimates above were compared with the source; this check does not certify the full proof. A direct invariant-mean argument can fail because a fixed partition level may have measure zero. The source instead chooses a finite averaging set A first and selects the level afterwards. Other OpenAI preprints are used only for corollaries.

\subsubsection{\texorpdfstring{287 \textperiodcentered{} All free group factors are isomorphic}{287 · All free group factors are isomorphic}}\label{fam:287}
\badges{Lean, 77 modules, in docs but not catalogue \badgesep{} combination \badgesep{} proves more than it states}

$L(\F_n)\cong L(\F_{n+1})$ for $n \ge 3$ by exhibiting a freely generating Haar n-tuple inside $L(\F_{n+1})$; $L(\F_2)\cong L(\F_3)$ via Dykema's amplification at $t=\sqrt{2}$; R\u{a}dulescu's dichotomy then makes all interpolated free group factors $L(\F_r)$, $1<r\le\infty$, isomorphic with fundamental group $\mathbb{R}_{>0}$. One consequence is that free entropy dimensions would not be generator invariants. The claim differs from expectations recorded in the initial extraction; Popa 2024 and Schafhauser--Tikuisis--White 2025 list the problem as open.

\paragraph{Extracted instruments.}

\begin{itemize}
\item \textbf{Small cocycle with a prescribed word value} (Lemma 4.1). Let $\Gamma=\F_n$, and place the coefficient h on the first generator only. The Fox cocycle of $w_m=(ab_1)(ab_2)\cdots(ab_m)a$, where $b_j=x_2^jx_3x_2^{-j}$, is $T_mh$ with $T_m=1+\sum\lambda(p_j)$ and freely generating $p_j$. Since $T_mT_m^*/m$ converges to free Poisson(1), which has no atom at 0, a cutoff inverse supplies h with $\|h\|_2\le(dm)^{-1/2}$ before finite-support approximation, where $d>0$ is a fixed spectral cutoff, and small $\|T_mh-\delta_e\|_2$. The spectral cutoff controls the discarded mass; the absence of an atom at zero does not supply a uniformly bounded inverse. The associated flow changes each generator by $O(\|h\|)$ while changing the long word by a full unitary factor. Fox calculus and the free Poisson law are classical; the initial extraction identified their use to obtain this displacement as the new feature. The repository\textquoteright{}s Lean development uses a combinatorial partial inverse and also treats $n=2$, beyond the paper\textquoteright{}s argument.
\item \textbf{Trace-preserving polynomial flow driven by conjugates of a bounded logarithm} (Prop. 3.1): $S=\arg C$, ODE $A_j'=i(\sum_gh_j(g)A_gSA_g^*)A_j$; trace preservation via a universal initial trace identity on a formal algebra; displacement bounds via a three-block free-Khintchine inequality. A specialization of Voiculescu's cyclomorphy to Haar unitaries with relations; incremental. The new choice is the fixed bounded logarithm and that the flow fixes C, so a Nielsen move converts the word's displacement into absorption of the extra generator.
\item \textbf{Adaptive error budgets and $L^2$-density $\Rightarrow W^*$-generation} (\S{}5): The current-stage estimate treats the newly produced word as one unitary slot. After substitution, the word-polynomial Lipschitz constant does depend on word length: $L(F)=\sum_\nu|c_\nu|\,|v_\nu|$ for $F=\sum_\nu c_\nu v_\nu$. Future error budgets are chosen after that constant is known; norm convergence preserves the Haar law. Incremental; Shlyakhtenko's Fuchsian preprint of 10 Sep 2026 (``obtained using ChatGPT Pro'') already uses finite automorphisms fixing one element plus limiting generation, and the Lean uses Baire instead.
\end{itemize}

\paragraph{Reused.}

Voiculescu's free probability; Fox calculus; Ricard--Xu Khintchine; cyclomorphy as a model; Nielsen moves; Dykema amplification and R\u{a}dulescu dichotomy (the Lean avoids both and proves $L(\F_3)\cong L(\F_\infty)$ directly).

\paragraph{Caveats.}

The claimed result would refute the Connes--Shlyakhtenko $L^2$-Betti conjecture for $L(\F_n)$; this consequence is not discussed in the source. A broader free-product absorption statement would require separate hypotheses and is not asserted here. The reasoning summary identifies trace preservation and limiting generation as dependencies.

\section{\texorpdfstring{Family by family, second group}{Family by family, second group}}\label{sec:cards2}

The second group belongs to the initial extraction and uses the same entry format. An entry marked ``no fresh search'' was compared with directly retrieved primary sources but did not receive a separate literature search.

\subsection{Number theory and arithmetic geometry}

\subsubsection{\texorpdfstring{002 \textperiodcentered{} Full BSD formula from low Selmer corank}{002 · Full BSD formula from low Selmer corank}}\label{fam:002}
\badges{not formalized \badgesep{} new construction \badgesep{} three-paper chain plus an unpublished preprint}

If the q-power Selmer group of $E/\mathbb{Q}$ has corank 0 or 1 at some prime q, then rank equals analytic rank equals that corank, Sha is finite, and the exact BSD leading-term formula holds with no hypotheses on reduction, torsion, isogenies, CM or residual image. The main paper's own contribution is the exact valuation at every odd prime; rank and finiteness come from the companion ``Selmer converse'', the 2-adic part from the ``two-primary'' companion, and the auxiliary quadratic fields from family 006, which rests on Smith's 2025 preprint.

\paragraph{Extracted instruments.}

\begin{itemize}
\item \textbf{Tame variables, finite models and residual concentration} (\S{}\S{}3--4). The coefficient ring is $\mathbb Z_p[[t,u_1,\dots,u_m]]$ where t is the cyclotomic or anticyclotomic variable and each $u_j$ is a nonprincipal-ultrafilter limit of cyclic p-power characters at auxiliary primes $r\to1$ p-adically. Neková\v{r}-style Selmer complexes are replaced by bounded finite free cochain models whose matrices have ultrafilter limits, keeping cup products, duality and homotopies. ``Residual concentration'' adjoins tame directions until one-sided determinants are nonzero mod (p,t), using a first-order lift test with Chebotarev-chosen commutators. This is where Eisenstein and small-image primes are absorbed. This extends Taylor--Wiles--Kisin auxiliary primes and Scholze--Manning ultrapatching; the combination (ultralimits of tame characters as formal Iwasawa variables, later specialized to zero) had no precedent in the literature examined.
\item \textbf{Integral character-division test and codimension-two divisor removal} (Lemma 3.9, \S{}9): Weierstrass division in a distinguished tame coordinate makes an analytic-over-algebraic quotient integral; any divisor would force a Selmer jump along a codimension-one locus, excluded by adjoining directions until simultaneous jumps have codimension $\ge 2$. This replaces main-conjecture input by divisibility plus divisor avoidance. Combination.
\item \textbf{Theta implication with a retained nilpotent tame parameter} (companion \S{}\S{}4--9, main \S{}\S{}7--8): a holomorphic theta family on U(3,1) at both real places of $\mathbb Q(\sqrt{DD'})$ whose constant Fourier--Jacobi term is a toric period specializing via BDP's inverse differential operator to a truncated series over $k[t]/(t^b)$; cusp lifting mod $p^n$, base change to $\mathrm{GL}_{4}$, Bellaïche--Chenevier extension extraction giving Artin length $\ge b$ against a first-order upper bound $\le 2$. This replaces Skinner--Urban and Castella--Liu--Wan Eisenstein congruences, with no ordinarity of E. New as a construction.
\end{itemize}

\paragraph{Reused.}

Gross--Zagier and Kolyvagin; Kato's zeta elements; Kolyvagin derivatives with Mazur--Rubin switches; BDP measures; Cassels and Milne isogeny invariance; Burungale--Flach for CM rank 0; Chai--Hida rigidity; Poitou--Tate.

\paragraph{Caveats.}

Proposition 5.14 asserts Kato's divisibility at every odd prime including Eisenstein primes, stronger than the statements found in the initial comparison, and is a principal step to check. The foundational steps require attention (ultrafilter limits of cochains, ``no $\mathrm{lim}^{1}$ term''). The companion's appendix cites an undated 185-page author-hosted draft, apparently not load-bearing.

\subsubsection{\texorpdfstring{007 \textperiodcentered{} Ordinary two-point Chowla}{007 · Ordinary two-point Chowla}}\label{fam:007}
\badges{repository reports Lean coverage; comparator and solution coverage not independently checked here \badgesep{} combination bordering new}

For fixed integers $a_1,a_2\ge1$ and $b_1,b_2\ge0$ satisfying $a_1b_2-a_2b_1\ne0$, the source claims
\[
 \left|\sum_{1\le n\le X}\lambda(a_1n+b_1)\lambda(a_2n+b_2)\right|
 \le C_{a_1,a_2,b_1,b_2}\frac{X}{(\log X)^c}\qquad(X\ge3).
\]
The exponent $c>0$ is absolute; the implied constant can depend on the fixed forms, and no uniformity for coefficients growing with $X$ is asserted. The corrected binary Elliott conclusion is qualitative, under uniform non-pretentiousness of an original factor and the same conditions on the affine forms. The initial extraction describes ordinary two-point Chowla as open before this release. The exponent is not given numerically.

\paragraph{Extracted instruments.}

\begin{itemize}
\item \textbf{Terminal-scale divisor bins with padding anti-concentration.} Dilation $n\mapsto qdn$ changes the scale, which explains the role of dilation-invariant logarithmic averaging in the earlier approaches of Tao, Helfgott--Radziwi\l{}\l{} and Pilatte. Here divisor pairs are confined to bins $j\eta\le\log(dq)\lt(j+1)\eta$ so the full-divisibility term is F(X)+O($\eta$) at a fixed scale; d is a product of J primes from disjoint supplies, q a squarefree ``padding'' weighted $4^{\omega(q)}$; a characteristic-function anti-concentration lemma spreads divisor mass so each bin costs $(C/\sqrt W)^J\le e^{-J}=(\log X)^{-c}$. This replaces dilation invariance of logarithmic weights. The fixed-scale bin plus mass-spreading padding is not in Helfgott--Radziwi\l{}\l{}, Pilatte or Tao--Teräväinen.
\item \textbf{Finite residue comparison via Braverman's theorem.} Prime sets have joint period far exceeding X, so CRT cannot compare interval averages with the independent-residue model. Each residue is encoded as the first $L^{15}$ bits of a jittered quotient, equality tests become DNFs, the bit law is made exactly t-wise uniform by the Alon--Goldreich--Mansour Fourier trick, and Braverman's theorem (polylog independence fools $\mathrm{AC}^{0}$) gives error $e^{-L^{10}}$ for every event in the proof, all written as bounded-depth circuits. This replaces Pilatte's rank-truncated sieve. The initial extraction found no earlier arithmetic application of Braverman's theorem.
\item \textbf{Degree-normalized two-band divisor graph with matched centering} (Part II): a core band weighted $e^{2\kappa}$ per prime and a center band with factor $1_{p|x}-\theta/p$, $\theta=1/4$ chosen so that lit and unlit residues cancel to $O(p^{-2})$ before absolute values; edges divided by the origin's total divisor weight. Combination over Helfgott--Radziwi\l{}\l{} and Pilatte's graphs.
\end{itemize}

\paragraph{Reused.}

Helfgott--Radziwi\l{}\l{} divisibility graph and prohibited paths; Pilatte's composite labels and non-backtracking operator; Matomäki--Radziwi\l{}\l{}--Tao short sums (corrected 2022 version); Granville--Soundararajan distance; Brun truncation.

\paragraph{Caveats.}

The reasoning summary makes the conclusion conditional on simultaneous-witness elimination, forest counting, finite-residue comparison and spectral-transfer bounds, while the paper states these steps as proved. They are the main points identified for further checking. Part II gives no rate.

\subsubsection{\texorpdfstring{010 \textperiodcentered{} Fontaine--Mazur modularity at the prime 2}{010 · Fontaine–Mazur modularity at the prime 2}}\label{fam:010}
\badges{not formalized \badgesep{} combination \badgesep{} input to Hilbert's tenth}

Every continuous irreducible representation $r:G_{\mathbb Q}\to\mathrm{GL}_2(\overline{\mathbb Q}_2)$ that is odd ($\det r(c)=-1$ for complex conjugation $c$), unramified outside finitely many primes, and de Rham at 2 with two distinct Hodge--Tate weights is a Tate twist of the 2-adic Galois representation attached to a classical cuspidal eigenform, ``including all residual representations'': scalar, reducible, dihedral. The earlier results cited by the initial extraction covered $p=2$ only for irreducible non-solvable and nearly-ordinary solvable residual images. Family 004 rests on this paper alone (not on its two satellites).

\paragraph{Extracted instruments.}

\begin{itemize}
\item \textbf{Marked lines and transverse charts} (\S{}\S{}5--6). After a solvable totally real base change through abelian 2-extensions with 2 split and a bounded number of bad primes, record an invariant ordinary line only at the $m=[F:A]$ dyadic places over one place of A, with $m/d\to0$. Restricting to points where marked lines are transverse to every global invariant line forces the local residual representation to be scalar at marks and makes the extension cocycle vanish there, so the reducible locus has dimension below the connectedness dimension (SGA2) of the chart and the residual fiber is connected. This replaces Skinner--Wiles and Pan's ordinary covering, which records a line at all p-adic places. The proposed new step is to mark a vanishing fraction of places to obtain connectedness, while leaving the propagation step unrestricted at 2.
\item \textbf{Finite-jet ultrapatching at a characteristic-two curve} (\S{}4): Taylor--Wiles patching at a one-dimensional prime with residue field k((t)), enhanced coefficient ring carried through all rings, ultraproducts of jets, a uniform Artin--Rees bound; full support obtained on the tame ``$\eta$ = i'' model and transferred to the unipotent ``$\eta$ = 1'' model via the common special fibre (Khare--Wintenberger's characteristic-2 form of Ihara avoidance). Incremental over Pan 2022 \S{}4; new only in the jet bookkeeping that survives a non-finite residue field.
\item \textbf{Ordinary character family with Thorne points and a polynomial-orbit lemma} (\S{}7): for locally reducible targets, an ordinary family finite dominant over weight space whose dense weight-$\{0,1\}$ points are modular by Thorne's Theorem D; a direct Banach-representation lemma shows a unit-eigenvector fixed by the unipotent radical is locally algebraic ``without assuming that algebraic vectors commute with specialization''. Combination; the seed set is entirely dependent on Thorne.
\end{itemize}

\paragraph{Reused.}

Chenevier determinants; Wang-Erickson algebraization; Böckle--Iyengar--Pa\v{s}k\={u}nas local geometry; Pa\v{s}k\={u}nas--Tung blocks; Colmez--Dospinescu--Pa\v{s}k\={u}nas; Breuil--Emerton; Newton--Thorne adjoint Selmer vanishing; Billerey--Menares lifts; Jacquet--Langlands and cyclic base change; Ramakrishna switching.

\paragraph{Caveats.}

The principal external input is Thorne's August 2026 preprint~\cite{src:thorne}; without that input, the stated argument for locally reducible targets, and the Hilbert's-tenth reduction that uses it, would require a replacement. Tabaâ 2020 (unrefereed) is cited for normality of completed tensor products. The pro-modularity satellite is doubly dependent. The paper consistently attributes each device to Pan, Khare--Wintenberger, Skinner--Wiles, BIP and Thorne; the new mathematics is combination and $p=2$ bookkeeping. The initial extraction completed fresh searches and reports checking the statements of each essential citation against its primary source.

\subsubsection{\texorpdfstring{040 \textperiodcentered{} Bloch's conjecture for surfaces with $p_g=0$}{040 · Bloch's conjecture for surfaces with p\_g = 0}}\label{fam:040}
\badges{not formalized \badgesep{} new strategy}

For every smooth connected projective complex surface with $p_g=q=0$, degree is an isomorphism $CH_0(S)\to\mathbb Z$ with integral coefficients; hence the full Bloch conjecture for $p_g=0$ (general type forces $q=0$, $\kappa<2$ is Bloch--Kas--Lieberman). For $p_g=0$, the stated motivic consequences concern the rational Chow motive $h(S)_{\mathbb Q}$: Kimura-finiteness and vanishing of its transcendental degree-two summand. The pure-Tate conclusions for $h(S)_{\mathbb Q}$, powers $S^m$ and Hilbert schemes $S^{[n]}$ additionally require $q=0$. The integral zero-cycle statement does not imply an integral Tate-motive decomposition. Two unrefereed competing preprints (Guletskii, Banerjee) are cited as ``proposed proofs'' with no logical dependence.

\paragraph{Extracted instruments.}

\begin{itemize}
\item \textbf{Marked Chow relation from Quot-scheme localization} (Thm chow-relation). On the Quot scheme of rank-0 quotients of $L_1\oplus L_2$ with fixed determinant, $\mathbb G_m$ acts with weights $\pm1$. The class $\Xi(t)=q_*(h^p\,ch_2(Q)_{(1)}ch_2(Q)_{(2)}\cap[B]^{vir}\times[X\times X])$ has equivariant codimension $p+4-v$; if $v>p+2$ its codimension-2 part would carry a negative power of t, impossible for a non-localized class, so it vanishes. Graber--Pandharipande localization rewrites the same class as a sum over fixed loci $|T|\times|H|\times X^{[n_1]}\times X^{[n_2]}$, yielding $\sum_WC_W=0$ in $CH^2(X\times X)_{\mathbb Q}$ with each $C_W=c_W[\Delta_X]+$ external products. This supplies Bloch--Srinivas decompositions needing no family, covering or automorphism. Positive-dimensional localization relations are standard (Faber--Pandharipande; Marian--Oprea--Pandharipande used this exact setup for Lehn's conjecture); the initial search found no earlier use of this approach for Bloch's conjecture.
\item \textbf{Chow universality with two output factors retained and the universal diagonal coefficient}: run the Ellingsrud--Göttsche--Lehn recursion in Chow without integrating two marked copies of X; components containing both outputs give a scalar times $[\Delta_X]$, computed by a formal rule depending only on $\int c_2$ and intersection numbers, hence transportable between surfaces with equal numerics ``whether or not T, H are effective''. The initial assessment is a new combination, with related work by Maulik--Negu\c{t} and Negu\c{t}--Oberdieck--Yin.
\item \textbf{Auxiliary-surface comparison via the divisor cosection} ($K^{2}=9$): on the Cartwright--Steger surface ($p_g=q=1$, same Chern numbers as a fake projective plane), Chang--Kiem's cosection eliminates every divisor component except the empty divisor and the unique canonical divisor, pairing with Re $\omega$ gives $c_{out}=s\,c_{mid}$, a length-2 blow-up computation gives $c_{mid}=12$, so $c=12(s+\tau)>0$, bypassing a length-20 Hilbert-scheme integral. The construction distinguishes $p_g$ exactly through whether $|K|$ is empty or a point. New as a strategy.
\end{itemize}

\paragraph{Reused.}

Bloch--Srinivas; Roitman torsion; Graber--Pandharipande; MOP/Oprea--Pandharipande Quot obstruction theories; EGL recursion; Carlsson--Okounkov; Toën--Vaquié and Schürg--Toën--Vezzosi derived moduli; Kiem--Li cosection; Serre's classification of indefinite unimodular lattices for $K^{2} \le 8$.

\paragraph{Caveats.}

The diagonal coefficient is a formal choice in a non-unique decomposition; detection is legitimate only because the test class is orthogonal to every divisor appearing, which the paper checks. Derived-geometric foundations (fixed obstruction theory, Manolache compatibility) are asserted in a few pages each. No fresh literature search was completed for two of the instrument comparisons.

\subsubsection{\texorpdfstring{193 \textperiodcentered{} Serre's positivity conjecture}{193 · Serre's positivity conjecture}}\label{fam:193}
\badges{not formalized \badgesep{} new device}

For a commutative Noetherian regular local ring $R$ of dimension $d$ and nonzero finitely generated $R$-modules $M,N$ with $\ell_R(M\otimes_R N)<\infty$ and $\dim_R M+\dim_R N=d$, the source claims $\chi^R(M,N)=\sum_{i=0}^d(-1)^i\ell_R(\mathrm{Tor}_i^R(M,N))>0$, without a restriction on characteristic or ramification. The previously open case was ramified mixed characteristic $V[[x]]/(p-f)$, $p\in\mathfrak m^2$, which Bhatt--Hochster--Ma still called ``the most tantalizing problem that remains''. Load-bearing external inputs (Cai--Lee--Ma--Schwede--Tucker perfectoidization results, Land--Tamme cdh descent for truncating invariants, Ma's Frobenius vanishing) were compared with their arXiv statements in the initial extraction.

\paragraph{Extracted instruments.}

\begin{itemize}
\item \textbf{Supported rational descent} (Thm 3.1): let $A$ be a commutative Noetherian $\mathbb Z_{(p)}$-algebra of finite Krull dimension, $J\subset A$ an ideal containing $p$, and $H$ a connective associative $A$-algebra with central $A$-action. Let $T$ be a finite commutative $A$-algebra with a finite group $G$ acting over $A$, such that $\operatorname{Spec}T\to\operatorname{Spec}A$ is surjective and $G$ is transitive on each geometric fiber. The claimed equivalence is $K^J(H)_{\mathbb Q}\simeq(K^J(H\otimes_A^{\mathbf L}T)_{\mathbb Q})^{hG}$. In the proof, the fiber of rationalized K-theory under inverting p is truncating because Goodwillie's HN terms cancel, hence a cdh sheaf; finite-flat descent by rational transfers; flattening by blow-up; a ``geometric orbit argument'' computes invariants. Combination: Voevodsky's h-topology and Suslin--Voevodsky's pseudo-Galois covering lemma (uncited) transported to spectrum-valued truncating invariants with noncommutative coefficients.
\item \textbf{Rank quotient and supported asymptotic complexes} (\S{}4): fix $s_n\to\infty$ and an ultrafilter; the perfect hull of sequences $(S^{a_n})$, $a_n=O(s_n)$, with arbitrary matrix sequences as morphisms, modulo the thick subcategory of $o(s_n)$ sequences. Homology functors land in Serre quotients by null-ultralimit normalized length, giving nonnegative Betti numbers and Euler maps; a Drinfeld dg quotient exhibits the category as Perf of a connective algebra so instrument 1 applies. This turns the sequence of truncated minimal resolutions of finite approximations to the perfectoid parameter quotient into one K-theory class. New as a device; relatives are Bhatt--Hochster--Ma's lim Cohen--Macaulay sequences, Schoutens's ultraproducts, Chuang--Lazarev rank functions.
\item \textbf{Two norm comparisons instead of invariance of individual evaluations}: ordinary and normalized perfectoid Euler characteristics are compared only after summing over Galois conjugates; both equal $r\chi_A$ on classes pulled back from the regular base, and descent says invariant rational classes come from A. ``No invariance of $\lambda$ under G is asserted or needed.'' This is the idea that bypasses non-functoriality of perfectoidization. Combination.
\end{itemize}

\paragraph{Reused.}

Serre's reduction to prime pairs via Roberts and Gillet--Soulé vanishing; Cohen--Gabber normalization; Faltings normalized length; CLMST perfectoidization in mixed characteristic and perfection plus Monsky and Roberts in characteristic p; Bhatt--Hochster--Ma's ``only $H_{0}$ survives'' mechanism via a Koszul-polynomial trick.

\paragraph{Caveats.}

The argument relies on several results from 2022--2026. Further checking should focus on the Verdier-quotient identification $\simeq\mathrm{Perf}(H)$ and its base change, spectral sequences on idempotent summands that need not lift to chain maps, and the interaction between ultralimits and ordinary limits. A December 2025 AI-generated note had already asserted the prime-quotient form; the paper treats it as a claim rather than an established theorem. The initial extraction retrieved primary sources directly but did not complete fresh searches for this family.

\subsection{Analysis, dynamics and mathematical physics}

\subsubsection{\texorpdfstring{143 \textperiodcentered{} Hilbert's sixteenth problem, uniform bounds}{143 · Hilbert's sixteenth problem, uniform bounds}}\label{fam:143}
\badges{not formalized (companion only) \badgesep{} combination}

For each degree d there is a finite B(d) bounding the number of limit cycles of every real planar polynomial vector field of degree $\le d$: the uniform-boundedness assertion in the second part of Hilbert's sixteenth problem. The bound is ineffective; Écalle--Ilyashenko individual finiteness is not used as a premise. The manuscript has 160 pages. The companion concerns the classical Liénard system $\dot x=y-F(x)$, $\dot y=-x$, with $\deg F\le5$; it claims at most two isolated periodic orbits. The repository supplies a formalization of this self-contained companion.

\paragraph{Extracted instruments.}

\begin{itemize}
\item \textbf{Packets on saturated logarithmic flags and the separation theorem} (Defs 2.1, 3.1, Thm 3.2): a flag $Z_1\gg\dots\gg Z_m\gg1$ of scales with exact triangular log relations; a packet is an analytic function plus two lateral trees of successive expansions in $e^{-aZ_j}$ on nested complex bands with transition estimates and antiholomorphic sources where determinations are glued; if all paired terminal coefficient germs vanish the function is identically zero, by weighted Phragmén--Lindelöf and Cauchy corrections. This is a multiscale, parametric, two-determination reorganization of Ilyashenko's superexact series and functional cochains, with ``differentiated closure'' and ``retestability'' added. Combination.
\item \textbf{Absolute isolated-zero finiteness and projection count} (Def 11.1, Thms 5.4, 11.2): an open analytic system with margins and its finite differential-algebraic closure has finitely many isolated real solutions when all coordinates including parameters are unknowns; with row rank this gives a uniform bound on fiber components (proper barrier, one generic tilt, Lagrange critical-point induction). Along an escaping sequence of isolated zeros, build a flag, eliminate with a deferred floor, apply separation, perturb a free scale to produce a non-isolated solution. This converts sequence-dependent finiteness into a coefficient-uniform bound. Incremental for the projection principle (Khovanskii 1984, Gabrielov); a new combination in its connection with packet quasianalyticity.
\item \textbf{Finite passage templates and matching systems} (Thms 9.5, 11.3): Lion--Rolin preparation of the slope field with all coefficients as parameters, in Euler clocks, puts every cell in an exact monomial normal form; Khovanskii's planar Rolle bounds boundary crossings so every cycle is a bounded word in a finite alphabet, matched end to end in finitely many analytic systems; hyperbolicity gives full rank and one component per cycle. This bypasses desingularization, polycycle classification and parameter compactness. Combination; the organizing idea of parametric preparation as a substitute for blow-up was not found elsewhere.
\end{itemize}

\paragraph{Reused.}

Lion--Rolin and Opris preparation; Bierstone--Milman and Gabrielov uniform finiteness; Khovanskii's Rolle; Gevrey truncation; Ilyashenko--Yakovenko cyclic systems; rotated vector fields to make cycles hyperbolic.

\paragraph{Caveats.}

The paper cites Yeung 2025 for a gap in Ilyashenko 1991 at a differential-closure step and claims to repair it through differentiated implicit calculus. That assertion and Sections 6--8, which devote 80 pages to checking packet hypotheses for four transfer models, merit detailed examination. Degenerate singularities and polycycles are not named outside the introduction. The literature comparison was limited.

\subsubsection{\texorpdfstring{263 \textperiodcentered{} The ionization conjecture}{263 · The ionization conjecture}}\label{fam:263}
\badges{headline not formalized; companions yes \badgesep{} new}

The source uses the nonrelativistic fermionic Coulomb Hamiltonian on $\bigwedge^n L^2(\mathbb R^3;\mathbb C^2)$, with two electron spin states and kinetic term $-\Delta/2$. The $M\ge1$ fixed nuclei have distinct positions $R_j$ and real charges $z_j\ge1$; $Z=\sum_{j=1}^M z_j$ is total nuclear charge. The system's net charge in electron-charge units is $Z-n$, and excess electron number is $n-Z$. Nuclear repulsion is omitted as a constant independent of $n$. For the ground-energy infima $E_n$ (with $E_0=0$), strict binding $E_n<E_{n-1}$ is claimed to imply $n\le Z+CM$, with universal $C$; threshold bound states are not excluded. For a neutral atom with one nucleus at the origin and integer $Z\ge1$, let $\rho_{\Psi_Z}$ be the spin-summed density of any normalized ground state, of total mass $Z$. Its outer half-electron radius is
\[
 R(\Psi_Z)=\inf\left\{r\ge0:\int_{|x|>r}\rho_{\Psi_Z}(x)\,dx\le\frac12\right\}.
\]
The stated bounds are $0<c\le R(\Psi_Z)\le C$ and $0<c\le I_1(Z):=E_Z(Z-1)-E_Z(Z)\le C$, uniformly in such $Z$. The radius retains a fixed exterior mass, not a fixed fraction of all electrons. The initial extraction lists repository formalizations for two Thomas--Fermi asymptotic companions (607 and 428 files), not an independently checked proof of the main strict-binding assertion.

\paragraph{Extracted instruments.}

\begin{itemize}
\item \textbf{Noisy-observation conditioning with a dimension-free tilt cost} (Lemmas 5.1--5.2). Observe every electron position with independent smooth noise of width $\ell_j=r_j^{1.01}$ at dyadic scales, forget labels, obtaining decreasing $\sigma$-fields. For a normalized eigenstate $\Psi$ and an observation event $A\in\mathcal F_j$ of probability $p>0$, put $p_A(x)=\Pr(A\mid x)$. With the stated noise law and $\sum_{h\ge j}\ell_h^{-2}\le C_\ell\ell_j^{-2}$, the vector $\Psi_A=\sqrt{p_A/p}\,\Psi$ is antisymmetric, in the form domain, has the conditional law, and $q[\Psi_A]\le e_\Psi+C\ell_j^{-2}\log^5(e/p)$, with no factor of n or of the number of arrays. The proof uses a score bound, Rademacher symmetrization and the eigenstate multiplier identity. Any failure of a local Thomas--Fermi comparison becomes a legitimate low-energy state. This replaces the deterministic Hartree--Fock density unavailable for correlated states. The initial comparison found no matching construction in the examined work of Lieb, LSST, Fefferman--Seco, SSS, Solovej, Nam or Lewin.
\item \textbf{Conditional density determines the field} (Lemmas 3.4--3.6, 5.3, Prop 10.1): a spatial cut leaves an antisymmetric conditional core with exact energy identity; coherent-state Bernoulli-occupied Slater trials (upper) and per-configuration Neumann bracketing through the empirical fine density of removed electrons (lower) share the same Thomas--Fermi constant; two posterior laws are compared via separate tower identities (``no equality of posterior laws is used''). Combination.
\item \textbf{Outward propagation by convexity-preserved subsolutions under observation forgetting} (Lemma 5.4, Prop 11.2): subsolution fields with reaction $4\pi kt_+^{3/2}$; Jensen for convex f and the tower identity preserve the inequality as observations are forgotten, errors paid by a density with summable expectation. The construction applies Solovej's multiscale approach to a martingale of conditional fields. Combination.
\end{itemize}

\paragraph{Reused.}

Lieb's 1984 bound; Zhislin binding; IMS localization; Lewin's sector localization; Lieb--Thirring; Lieb--Simon Neumann bracketing and Thomas--Fermi theory; Solovej's coherent states and exterior screened-potential strategy; Bénilan--Brezis.

\paragraph{Caveats.}

The argument has many quantitative parameters ($L=10^6$, exponents 1.01, probabilities $r^{25}$, $r^{40}$); the full bookkeeping was not checked. The formalization catalogue marks ``Partial progress'', review ``unchecked'', and omits the radii development from its results list. The literature comparison was limited.

\subsubsection{\texorpdfstring{265 \textperiodcentered{} Two-dimensional area law from a global gap}{265 · Two-dimensional area law from a global gap}}\label{fam:265}
\badges{not formalized \badgesep{} new}

Let $\Lambda$ be a finite induced square-lattice graph, with local dimension $q\ge1$, and let $H=\sum_Xh_X$ have one Hermitian term per nonempty support $X$, with $\|h_X\|\le J$ and support diameter at most $R$ in the graph distance of $\Lambda$. Assume a unique ground state and a full-system gap $\Delta>0$. The source claims $S(A)\le C(q,R,J,\Delta)|\partial_\Lambda A|$ for every $A\subseteq\Lambda$, where $\partial_\Lambda A$ is the set of lattice edges crossing from $A$ to its complement within $\Lambda$. No gap on restrictions, frustration-freeness or translation invariance is assumed. For the open $L\times L$ square grid with bounded on-site and nearest-neighbor interactions, a companion gives a normalized PEPS approximation with vector error at most $L^{-1}$ and bond dimension at most $CL^c$, with $C,c$ depending on $q,J,\Delta$. The constant is existential. Earlier proofs cited in the initial extraction (Anshu--Arad--Gosset 2022) need frustration-freeness.

\paragraph{Extracted instruments.}

\begin{itemize}
\item \textbf{Optimized nested positive filters with stationarity and Schmidt pinning} (\S{}3). Maximize $\|L_m\cdots L_1\Omega\|$ over positive filters on nested regions with constraints $\mathrm{Tr}L_j^{2/a_j}=1$ and harmonic weights; stationarity forces each filter to commute with the filtered marginal, so filters are clipped marginal powers; a cosh estimate bounds excitation energy, and the global gap converts this into an entropy sum (via a tilt by $\rho^{-u/2}$ and the three-lines theorem). The opposite bound pins a Schmidt vector in trial filters with analytic continuation. Subtracting gives a conditional-entropy discount and a subvolume exponent with no local gap. This replaces the AGSP subvolume argument and is classified as new in the initial extraction.
\item \textbf{Replica metric transport} (\S{}\S{}5--8): on $\mathrm{Sym}^k(V)$, central operators $W_Q\approx e^{-tF_Q}$ with Schur labels $\approx k$\textperiodcentered{}S(Q), equal on complementary regions; band metrics combined through finite trees of Kubo--Ando geometric means as a small subsystem switches sides. The estimates are $-\partial_p\log N^2\ge k\,a\,\langle I(x{:}F|P)\rangle-$err, and energy $\le 2E_0+Ca^2\cdot(\text{split probability})$, so the displayed information term is linear in a, while the additional energy cost is quadratic in a and linear in the split probability; a slow scan through a collar leaves sublinear mutual information. This replaces Chebyshev/AGSP Schmidt-rank bookkeeping. In the literature examined, encoding conditional mutual information as the derivative of a log-norm comparable to energy had no precedent.
\item \textbf{Positive replacement, randomized scanner, amplification} (\S{}\S{}4, 9, 10): absolute values of Kitaev's filtered terms retain the gap with exponential tails; random schedules make each interaction rarely split; Poisson-clock products give a local unitary witnessing decoupling. The initial extraction classifies this as a new combination; taking absolute values is an elementary but essential step.
\end{itemize}

\paragraph{Reused.}

Lieb--Robinson bounds; Kitaev filter; Alicki--Fannes--Winter; strong subadditivity; Schur--Weyl; de Finetti; Kubo--Ando means; Schatten interpolation; random-order chain rule.

\paragraph{Caveats.}

The exponent estimates have a small margin ($\varepsilon=10^{-5}$, 0.999952 versus 0.99999), and the replica errors must be controlled uniformly over all histories. The argument gives no quantitative $\Delta$-dependence. Not every instrument received a fresh literature search.

\subsubsection{\texorpdfstring{268 \textperiodcentered{} The spin-one Haldane gap}{268 · The spin-one Haldane gap}}\label{fam:268}
\badges{not formalized; computer-assisted \badgesep{} new criterion}

For the pure spin-1 Heisenberg antiferromagnet on even periodic rings, the source claims a unique ground state for all even $L \ge 60$ and a uniform gap: $\gamma_L>\frac4{105}\log\frac{80}{79}>4.79\times10^{-4}$ for $L \ge 60$ and $>3.8\times10^{-3}$ for $L \ge 2304$, versus the numerical 0.4105. Odd rings are not covered. A companion discharges Tasaki's Assumption 2 for boundary-field chains. Haldane's 1983 conjecture for the actual Heisenberg chain had no rigorous proof; the chain is not frustration-free, so the usual AKLT and frustration-free finite-size criteria do not directly establish this claim.

\paragraph{Extracted instruments.}

\begin{itemize}
\item \textbf{Coupled purity bootstrap} (Prop 4.4). Two purities of one partition function $Z_n(b)$ are the Gibbs purity $T(n,b)=Z_n(2b)/Z_n(b)^2$ and the spatial purity $S(n,b)=Z_{2n}(b)/Z_n(b)^2$, formed from the spatial-transfer eigenvalues. If a probability vector has purity $P$, its normalized square has purity $P'$ satisfying
\[
 1-P'\le\frac{(1-P)^2}{2P^2}.
\]
This gives quadratic improvement in a high-purity regime, not convergence to purity one from arbitrary initial data (a uniform probability vector is unchanged). Squaring Gibbs weights doubles $b$; squaring spatial weights doubles $n$. The exact coupling identities remain
\[
 S(n,2b)=\frac{S(n,b)^2T(2n,b)}{T(n,b)^2},\qquad
 T(4n,b)=\frac{T(2n,b)^2S(2n,2b)}{S(2n,b)^2}.
\]
The source's quantitative bootstrap assumes even $n_0\ge4$, $b_0>0$, $0<u_0<1/6$, and constants $C,r$ such that
\[
\begin{gathered}
 G(u_0)\le C,\quad r=Cu_0<1,\quad C(1-3u_0)>3,\\
 G(u)=\tfrac12\bigl((1-u)^{-1}+(1-u)^{-2}+(1-u)^{-3}\bigr)^2,\\
 S(n_0,b_0)\ge1-u_0,\qquad T(2n_0,b_0)\ge1-u_0.
\end{gathered}
\]
Together with the transfer representation, these hypotheses give dyadic defect decay and interpolation to intermediate even lengths. The inequality $T(L,b)>1/2$ alone forces ground-state uniqueness: two ground states would have equal maximal Gibbs weights. The uniform gap $\gamma_L>\log(1/r)/b_0$ uses the quantitative decay and temperature bounds as well. The criterion uses no local ground-state projections; its novelty judgment concerns this gap mechanism, not the elementary purity identity alone.
\item \textbf{Self-adjoint spatial transfer operator on bond histories} (Prop 3.1): the Aizenman--Nachtergaele Poisson expansion regrouped by bond gives a real symmetric Hilbert--Schmidt kernel; twisted traces split the spectrum into multiplicity sectors. This gives a rigorous continuous-time quantum transfer matrix in the tradition of Suzuki and Klümper. Combination.
\item \textbf{Moment-to-spectrum extraction with exact-integer certificates}: partition functions of $\le 12$ sites enclosed by floor-Horner in big integers; polynomial filters and capped masses convert moments into spectral bounds; an SU(2)-structured periodic MPS certifies the energy density. The initial extraction classifies this computational method as incremental.
\end{itemize}

\paragraph{Reused.}

Aizenman--Nachtergaele 1994; Suzuki's quantum transfer matrix; Fannes--Nachtergaele--Werner MPS; Fox--Glynn truncation; Tasaki's 2025 index theorem (companion).

\paragraph{Caveats.}

The argument uses 27 thermal cases and 2 MPS contractions. The verifier was not run for this survey. Self-adjointness requires the real Cartesian decomposition, which restricts the scope of the criterion.

\subsection{Combinatorics, geometry and complexity}

\subsubsection{\texorpdfstring{103 \textperiodcentered{} $\mathrm{L}=\mathrm{RL}=\mathrm{BPL}$}{103 · L = RL = BPL}}\label{fam:103}
\badges{Lean, multi-tape Turing-machine statement \badgesep{} new}

Deterministic logspace equals randomized logspace, one- and two-sided error. The comparator states exactly $\mathrm{L}=\mathrm{RL}\wedge\mathrm{RL}=\mathrm{BPL}$ for a multi-tape model with fresh coins per step and polynomial halting on every coin tape. The record stood at Hoza's $\log^{3/2}n/\sqrt{\log\log n}$. The formalization has a narrower scope than the full construction: the effective compiler and numerical bounds are not formalized, and Lean identifiers are obfuscated, which limits the textual inspection reported here to the statement.

\paragraph{Extracted instruments.}

\begin{itemize}
\item \textbf{Exact correction and copy hierarchy} (\S{}4). From the nilpotent substochastic transition S on the configuration graph, $p_0=(I-S)^{-1}e$. A stage applies a correction E with $I-D=(I+E)(I-C)$, moving transition mass into the reward; E is a Lipschitz cutoff iterated a constant number of times (a ``factorial truncation'' lemma); ``detectors'' set column multipliers to suppress correction on overloaded columns; a copy step splits large entries among copies with an intertwining that transports rewards. Each stage cuts the max entry by $2^H$ while at most doubling active vertices; after O(log n) stages all entries are small. This replaces repeated squaring (Saks--Zhou), Richardson iteration and PRG-fooled branching programs. Shortcutting and vertex splitting have precedents; the initial comparison found no matching hierarchy with this doubling bound.
\item \textbf{Rank channels and conditional averaging in a shared environment} (\S{}\S{}3, 5): all randomness is O(log n) bits (a few matrices over $\mathbb F_P$ plus salts); retention of a vertex at nested rates via an $h_{0}$-wise-independent polynomial sampler; averaging over a hash slice without storing identifiers by a lazy walk on the slice via the affine action of $SL_H(\mathbb Z)\ltimes\mathbb Z^H$, with Shalom's explicit Kazhdan constants giving a spectral gap uniform in moduli and ``validated walk variants'' guessing the wrap vector. Combination: Margulis--Gabber--Galil lifted to H dimensions; new is the conditional sampler on hash slices.
\item \textbf{Allowance/clock controller with catalytic key equality} (\S{}\S{}7--11): telescoping multilevel estimates with longer walks at larger slack; rows stored as exceptions to a default; Schwartz--Zippel fingerprints; contour traversal of a program tree charged so the stack telescopes to O(B); a shared B-bit catalytic vector lets equality of two keys be XOR-accumulated in four excursions with O(1) local state; errors measured by a root-Lipschitz loss so relative errors compose. Combination; the catalytic-equality gadget is a candidate new component.
\end{itemize}

\paragraph{Reused.}

Nilpotent inverse of the configuration graph; k-wise independent hashing; relative property (T) expanders (Shalom, Kassabov); Riesz--Thorin; Carter--Wegman fingerprints; Cook--McKenzie traversal; Lange--McKenzie--Tapp reversibility; catalytic computation (BCKLS).

\paragraph{Caveats.}

The source contains about 400 KB of text and an eight-step chain of sufficiently large constants. It avoids storing a seed hierarchy, the source of the $\log^{3/2}$ space bound in the comparison made by the initial extraction.

\subsubsection{\texorpdfstring{112 \textperiodcentered{} Depth-three circuits beyond the square root}{112 · Depth-three circuits beyond the square root}}\label{fam:112}
\badges{Lean, 179 files \badgesep{} new}

The source gives a fixed language in P such that for every real constant $A>0$, $\Sigma_{3}$ circuits (unbounded fan-in, no bottom fan-in restriction, size = gate count) computing its n-bit slice need more than $2^{A\sqrt n}$ gates for large n: the explicit $2^{\omega(\sqrt n)}$ problem of Håstad--Jukna--Pudlák 1995 and Gurumukhani et al. 2024, record $2^{1.282\sqrt n}$ (Paturi--Pudlák--Saks--Zane 2005). The growth rate is an inverse tower in log n; this is weaker than $2^{n^{1/2+\varepsilon}}$, as the paper says.

\paragraph{Extracted instruments.}

\begin{itemize}
\item \textbf{Restriction bound independent of clause count} (Lemma 3.1). A width-k CNF H under a p-random restriction has the expansion $H|_\rho=K+K\sum_j(-1)^j\sum_{\text{paths}}\prod v_{i_h}/\prod a_{U_h}$. Here K conjoins the initially easy clauses ($\le b$ live variables), each path records successive hard clauses, and the denominators count violated easy clauses. A marker-mixture identity expresses each path term as a nonnegative weight times the expectation of a width-b CNF: randomly order clauses and markers, then retain the clauses preceding the first marker. The resulting correlation estimate is $|E[G\cdot H|_\rho]|\le W\cdot\mathrm{corr}_b(G)$ with $E[W]\le1/(1-\theta)$. The Poisson-type tail $\theta$ is $<$ 1 for any constant pk when b $\gg$ pk. Classical Håstad switching has a useful decay estimate when $pk<1/5$; the scale $p\approx1/k$ with $k\approx\log S$ leads to the $\sqrt{n}$ exponent. In the signed expansion, normalization averages over eligible violated clause indices and cancels their count. The initial extraction classifies this mechanism as new. Tal\textquoteright{}s $L^{1}$ switching lemma and small-restriction switching lemmas have different conclusions.
\item \textbf{Sparse partition} (Lemma 4.1). Every width-b CNF on m variables is a pointwise sum of $\le 2^{\varepsilon m}$ pairwise disjoint width-b CNFs, each with $C_{b,\varepsilon}m$ clauses, for fixed integer $b\ge1$ and $\varepsilon>0$. This follows by combining Impagliazzo--Paturi--Zane sparsification with the disjointness construction of Dell--Husfeldt--Marx--Taslaman--Wahlén. The source describes this component as incremental.
\item \textbf{Slice language} (Def 5.1). The input contains data x, Toeplitz hash bits, a monic polynomial, and t field elements. Acceptance requires the degree-zero coefficient of the polynomial evaluated at the hash to vanish in the quotient-ring representation. In the proof, an irreducible modulus is selected and the hash is required to be injective on the restricted data cube. Random independent field coefficients then give $t$-wise independent signs on that cube. Averaging subsequently fixes one setting of all non-data bits with a small failure probability over restrictions; a fixed deterministic slice is not itself a random $t$-wise independent family. The construction combines Krawczyk hashing, Alon--Babai--Itai limited independence, and the evaluator approach of Paturi--Saks--Zane. On the good random restrictions with $m$ surviving data bits, its correlation with every width-$b$ CNF is at most $2^{-m/8}$; the exponent $1/8$ is independent of fixed $b$, while the sufficiently large length threshold may depend on $b$. Parity lacks this property, so the argument is consistent with the tight parity bound.
\end{itemize}

\paragraph{Reused.}

Håstad-style restrictions and wide-clause truncation; Razborov's reverse encoding; sunflower lemma; universal hashing; t-th moment bounds.

\paragraph{Caveats.}

The argument does not depend on other OpenAI manuscripts. The mathematical review checked selected identities in Lemma 3.1 and the arithmetic in \S{}6. A separate literature comparison for the marker-mixture device remains incomplete.

\subsubsection{\texorpdfstring{156 \textperiodcentered{} Borsuk's conjecture fails in dimension 9}{156 · Borsuk's conjecture fails in dimension 9}}\label{fam:156}
\badges{Lean, 275 files, exact real arithmetic \badgesep{} combination leaning new}

The set X of rank-one projectors $uu^{\mathsf T}$, $|u|=1$, in the trace-one affine hyperplane of $\mathrm{Sym}_4(\mathbb{R})$, which is isometric to $\mathbb{R}^9$ under the Frobenius metric, has diameter $\sqrt{2}$ and cannot be covered by ten sets of smaller diameter, so $b(X)\ge 11>10$; by adjoining points, Borsuk fails for all $d \ge 9$. The previous range cited by the initial extraction was $d \ge 64$ (Jenrich--Brouwer). The witness is a compact copy of $\mathbb{R}P^3$, not a finite set; no finite 9-dimensional counterexample follows.

\paragraph{Extracted instruments.}

\begin{itemize}
\item \textbf{Odd spectral extension of an admissible map.} From a smooth partition of unity f on $\mathbb{R}P^3$ whose positive supports are disjoint for orthogonal lines, define $H(v)=|v|^2f([v])$, $E(Q)=k\int H(Q^{1/2}u)\,d\mu(u)$ on PSD matrices, and $F(A)=E(A^+)-E(A^-)$; F is odd with $\|F(A)\|_1=\mathrm{tr}|A|$, so it is an odd map between spheres and Borsuk--Ulam gives $m \ge N_k$ immediately, exactly tight. This test map supports the rest of the argument and replaces Frankl--Wilson and strongly-regular-graph counting. Combination; Kalai's survey names the set but asks no topological question about it; Adams--Bush--Frick 2021 is the nearest relative.
\item \textbf{Equality-case degree analysis.} At $m=N_k$ the map has odd mod-2 degree; a maximal support with fewer than k labels gives a positive-dimensional level set inside one affine chart, where the extension descends to $L\times\mathbb RP^{h-1}\to\mathbb RP^{a-1}$ with vanishing Stiefel--Whitney power, degree 0, contradiction. So every facet has exactly k labels and the facets form an $\mathbb F_2$ cycle. This converts a tight Borsuk--Ulam bound into rigid combinatorics. New as a technique for Borsuk's problem, as far as could be checked (no dedicated literature search was completed). The Lean replaces this step by an equivariant simplicial Smith-theory ``characteristic degree'', a different formal route.
\item \textbf{Six-label triangle systems} (\S{}\S{}5--6): the forced structure on the complement of any tetrahedron is, by an enumeration reported in the initial extraction, exactly the hemi-icosahedron $\mathbb{R}P^2_6$ (unnamed in the paper); the 10-label contradiction is a Walkup-style nonexistence of a 10-vertex neighborly 3-pseudomanifold with all facet complements $\mathbb{R}P^2_6$. Incremental; the Lean proves it by 145 machine-generated case files.
\end{itemize}

\paragraph{Reused.}

Projector embedding (Conway--Hardin--Sloane); partitions of unity; Sard; Borsuk--Ulam; mod-2 local degree; Stiefel--Whitney classes; Hinrichs--Richter point adjoining.

\paragraph{Caveats.}

The counterexample is continuous and does not settle finite-set Borsuk in $\mathbb{R}^{9}$. Its convex hull gives a convex-body counterexample, a consequence not discussed in the source. The paper cites an unrefereed AI-assisted GitHub note by Grinsztajn concerning dimension 63, without using it as an input. The literature comparison was limited.

\subsubsection{\texorpdfstring{266 \textperiodcentered{} Three mutually unbiased bases in dimension six}{266 · Three mutually unbiased bases in dimension six}}\label{fam:266}
\badges{headline not formalized; $N(6) \le 5$ formalized \badgesep{} combination \badgesep{} arithmetic and compiler assumptions}

The source claims $N(6)=3$ (Zauner's conjecture), as a theorem explicitly labelled ``computer-assisted'' and conditional on arithmetic and compiler conditions: a C++ pipeline in binary64 with hand-proved one-sided rounding margins, 4,470 seconds on 4 Oct 2026. The headline computation uses neither interval arithmetic nor exact algebra. The companion proves $N(6) \le 5$ exactly, the Matolcsi--Ruzsa--Weiner Fourier-vanishing conjecture for order-6 Hadamard matrices outside Tao's $S_6$ equivalence class, and Phangara's conjecture, with an exact rational verifier; the Lean covers the companion.

\paragraph{Extracted instruments.}

\begin{itemize}
\item \textbf{Seed ball plus fixed-frame completion via a U(2) polar-factor cover} (\S{}\S{}2--3). Dephase H; cover the nine phases of rows and columns 0--1 by 249,878 balls of radius 0.032 after a minimum-circular-gap normalization; inside a ball write H in 2+4 blocks, the kernels of the $2\times4$ blocks are 2-dimensional, the compressed unknown $Y=K_C^*DK_B$ is $2\times2$ with $\mathrm{tr}(6I-Y^*Y)\lt0.32$, hence within p of $\sqrt{6}$ times a unitary $V \in U(2)$; an Euler-angle box cover of U(2) with uniform Taylor bounds yields necessary tests. Every order-6 complex Hadamard matrix is represented by 13 parameters with no family assumption, giving a finite cover of the unclassified 25-phase continuum. Jaming--Matolcsi--Móra proposed discretizing all 25 phases without completing that computation in the work compared by the initial extraction. Combination (polar decomposition, box covers, dual certificates).
\item \textbf{Root certificates and tangent/normal refinement} (\S{}\S{}4--5): a center, a numerical tangent frame and constants such that every exact solution displacement satisfies a normal-component bound, from a residual identity with an approximate left inverse of the Jacobian; ``no manifold parametrization is required''. Incremental validated numerics without interval arithmetic.
\item \textbf{Injective moments} (companion): $\lambda(E)=E[\mathrm{Re}\prod H(\sigma(i),\tau(j))^{E_{ij}}]$ over random injections; orthogonality contracts order-one rows; an exact identity $3|\tilde\gamma|^2+\mathrm{avg}|g(\alpha)g(\beta)|^2=0$ forces either all permuted Fourier values to vanish or $H\circ H$ Hadamard, Newton's identities force cube roots, and a five-cycle gives Tao's $S_{6}$. A Delsarte-type moment certificate in the Matolcsi--Weiner lineage; the Fourier-vanishing claim applies outside the exceptional equivalence class of Tao's $S_6$ (previously Karlsson's family only). Combination leaning new.
\end{itemize}

\paragraph{Reused.}

JMMSW phase discretization and clique search; Rump's verified norm bounds; Higham's rounding analysis; Klappenecker--Rötteler lower bound; Newton's identities; Zauner's 2-design identity; Weiner's completion theorem; Peyrl--Parrilo rational certificates.

\paragraph{Caveats.}

The headline rests on about 70 pages of hand-derived rounding margins across seven C++ files, not independently rerun or verified for this survey; one condition is checked by inspecting disassembly addresses. The 24 Sept manuscript's proof rests on a 4 Oct run; an earlier run's evidence packet ``was not recovered''. The August 2026 classification preprint of Cárdenes Wuttig--Tindall is cited but is not used in the proof.

\subsection{Operator algebras, topology and symplectic geometry}

\subsubsection{\texorpdfstring{221 \textperiodcentered{} Mézard--Parisi formula for diluted spin glasses}{221 · Mézard–Parisi formula for diluted spin glasses}}\label{fam:221}
\badges{Lean, 405 files, not in catalogue \badgesep{} new}

For Poisson-diluted even-arity models whose interaction law satisfies the Panchenko--Talagrand factorization and positivity hypotheses (model-class conditions, verified for Viana--Bray and soft K-SAT, not assumptions on the Gibbs measure), the free energy equals the infimum over finite-depth hierarchical trial laws, with only first-moment integrability. The lower bound is stated for every bounded interaction whose law is permutation invariant, without the even-arity or factorization restrictions needed for the matching upper bound. Under the stated model-class hypotheses, no additional Gibbs-measure factorization assumption is imposed; establishing that factorization is the step Panchenko's conditional results left open.

\paragraph{Extracted instruments.}

\begin{itemize}
\item \textbf{Marked-depth Poisson perturbations with prescribed-tree identities} (\S{}7, Thm 7.6). Add $N^{3/4}$ bounded single-site factors $1+t_\nu g\sigma_i+u_\nu g_0S_0$ where g, $g_{0}$ are $\pm1$ functions of Rademacher marks sampled at prescribed depths of an auxiliary L-level hierarchy; the partition function iterates by power means with $m_d=d/L$. Scores concentrate (Efron--Stein, as in Ghirlanda--Guerra); a Poisson add-one identity and a two-scale probe $t=1/j$, $u\in[j^{-j},2j^{-j}]$ force a limiting analytic function to vanish identically. Shared marks act as sharing indicators, and polarization selects any branching pattern with exact signed coefficients. This extends Ghirlanda--Guerra, Barra--De Sanctis stochastic stability, Barbier--Panchenko side observations and Mourrat--Panchenko's enriched hierarchy (the latter two named in the trace, uncited in the paper). Combination.
\item \textbf{Depth shift with charging fresh splits} (\S{}8, Lemma 4.5, Thm 8.2). To extract conditional covariances one divides by coefficients of order $\delta^{2b}$; naively the error sum is O(1). Shifting every target branching depth by one level off all old depths forces each target vertex to be created by a fresh unary split of coefficient exactly $-\delta$, so $\sum|c_E|\le C\delta^{2b}$ and the ratio is bounded. The cost is controlled only on average over regular depth assignments via a vector-martingale energy bound; a three-copy estimate gives induction on leaf count; the limits are strictly ordered $N\to\infty$, then $L\to\infty$, then $\eta\to0$. This proves that all finite multioverlaps concentrate given the auxiliary hierarchy, the ``factorization'' Panchenko 2016 called missing. New, as far as could be determined; the trace marks it as the turning point.
\item \textbf{Recursive law encoding of a single site} (Lemma 9.1): conditional magnetizations of a uniform site encode an admissible Panchenko--Talagrand trial law reproducing all prescribed-tree moments. Incremental, a finite-N form of Austin--Panchenko.
\end{itemize}

\paragraph{Reused.}

Guerra, Franz--Leone and Panchenko--Talagrand interpolation (upper bound, reproved); Aizenman--Sims--Starr cavity; Panchenko's spin-distribution viewpoint; Ruelle probability cascades; Efron--Stein.

\paragraph{Caveats.}

The architecture and selected estimates in \S{}8 were inspected, without checking every inequality. The repository supplies a Lean development with review status ``unchecked''; it was not compiled here. A textual search for admitted proofs would not establish kernel verification.

\subsubsection{\texorpdfstring{285 \textperiodcentered{} Baum--Connes and Kadison--Kaplansky counterexamples}{285 · Baum–Connes and Kadison–Kaplansky counterexamples}}\label{fam:285}
\badges{not formalized \badgesep{} combination \badgesep{} reduced assembly-map scope}

The source constructs a finitely generated torsion-free group $G_{inj}$, an injection $\mathbb{Z}^2\to G$ and a degree-2 class pairing to 1 with the torus, such that the Bott class dies in $K_0(C^*_rG)$: rational injectivity of the reduced, coefficient-free assembly map fails. The claim does not contradict classical Novikov or maximal injectivity (Hanke--Schick). The source asserts finite generation; it establishes neither finite presentation nor its failure. The corrected Antonini--Buss--Engel--Siebenand result concerns specified exotic group $C^*$-algebras for finitely presented groups; identifying its target with the reduced algebra needs an additional condition, such as exactness, which is also not supplied here. Companions: a torsion-free group with a projection of trace in $(0,1/2)$ (surjectivity failure, Kadison--Kaplansky fails) and a group with p-torsion and an irrational-trace projection.

\paragraph{Extracted instruments.}

\begin{itemize}
\item \textbf{Heisenberg voltage tower with simultaneously sampled Bott symbols} (\S{}\S{}1--4). The construction begins with a layered bipartite graph labelled by interleaved bit streams with deterministic recovery (Osajda-type word uniqueness with an alphabet small enough that $\sqrt{|S|}/q\to0$); permutation tables chosen by an LLL-type sequential law; pruning of ``coherent returns'' and ``bicycles''. Each edge gets a voltage in the integral Heisenberg group H and a matrix weight from one distribution whose mean realizes a prescribed torus symbol spread along the central direction, while diffuse horizontal coordinates give every short non-law nontrivial holonomy (Schwartz--Zippel on Heisenberg words). Every tower layer carries $\ell^2(H)$ and incidences become Bott-symbol-valued band operators; the cover is the graphical relator. This extends Gromov--Osajda labelled-expander relators (no group-valued coordinates) and Gross--Tucker voltage graphs (no analytics). Combination; no prior use of a nilpotent voltage cover of an expander tower as a relator to transport deck-group K-theory was found.
\item \textbf{Parity-swindle Bott cancellation in a norm-sequence quotient} (\S{}6): four bounded profile families realized by finite Cayley polynomials, two explicit projection paths giving equal $K_{0}$ classes, finite additivity forcing the height-zero class to vanish; passage through $\prod C^*_r/\bigoplus C^*_r$ and a ``finite witness'' lemma. The paper calls it the Eilenberg-swindle pattern of coarse operator K-theory (Higson--Pedersen--Roe). The underlying idea is incremental; its use in this setting is classified as a new combination.
\item \textbf{Sheet-cone complex, exact cone chains, prime-order obstruction, Shapiro extension} (\S{}5b): cone off sheets, finite peeling gives exactness, a prime-order subgroup would contradict periodicity, so $\Gamma$ is torsion-free; cochain-level Shapiro extends the $H^{2}$ generator. Incremental.
\end{itemize}

\paragraph{Reused.}

Baum--Connes--Higson; Lück's trace theorem; an integrality attribution to Bartels--Greenlees--Wegner left unresolved by the initial extraction; Hanke--Schick; Emerson--Hudson; Kasparov descent; graphical small cancellation; Osajda labelling; LLL; Rieffel projections; Dykema--Rørdam.

\paragraph{Caveats.}

The construction uses large parameters (s a multiple of 20,000, scales $2^{500H_i}$). Their order of choice and the claimed absence of circular dependence require close checking. All three papers share the architecture, so a defect in it could affect the whole family. The main paper adapts and reproves methods from the Kadison--Kaplansky companion, which cites an unrefereed GitHub note as motivation. The comparison used directly retrieved primary sources without a fresh literature search.

\subsubsection{\texorpdfstring{293 \textperiodcentered{} Hyperinvariant subspaces and transitive algebras}{293 · Hyperinvariant subspaces and transitive algebras}}\label{fam:293}
\badges{Lean, 1,108-module import closure \badgesep{} new}

On every infinite-dimensional separable complex Hilbert space there is a nonzero bounded quasinilpotent operator S ($\|S^n\|\le\exp(-10\lfloor(n+1)^2/4\rfloor)$) with no nontrivial closed hyperinvariant subspace; its commutant is a proper strongly closed transitive algebra, refuting Kadison's transitive-algebra problem as well, which the catalogue summary omits. Companion: for every irrational rotation, an operator in the hyperfinite $\mathrm{II}_{1}$ factor with no nontrivial invariant projection in the factor, answering Zhu--Fang--Shi and the point-mass Brown-measure case left open by Haagerup--Schultz. Nothing is said about the ordinary invariant subspace problem.

\paragraph{Extracted instruments.}

\begin{itemize}
\item \textbf{Backward intertwiner over an ergodic base.} $S=\int^\oplus S_x$ over the 2-adic integers with Haar measure and the odometer $x\mapsto x+1$; each fiber a unicellular bilateral weighted shift with quadratic log-weights plus residue corrections so the cocycle drift has opposite signs at the two ends. V sends $e_j$ at x to $b_j(x)e_{j-1}$ at x+1, where b solves a commuting-square recurrence, normalized per fiber by a supremum finite a.e. by Borel--Cantelli but not integrable. Multipliers $L^\infty(\mathbb Z_2)\subset\{S\}'$ force any hyperinvariant subspace to be a measurable field of tails $K_{\ge s(x)}$; $V\in\{S\}'$ gives $s(x+1)\le s(x)-1$, impossible on a probability space. This extends Atzmon--Sodin's completely indecomposable shifts by transport between fibers. New as a construction; the Bishop-operator literature works in one fiber on invariant, not hyperinvariant, subspaces.
\item \textbf{Endpoint-conditioned convolution obstruction} (companion, Lemmas 4.1, 4.4): from a putative projection p of trace in $(0,1)$, exact convolution identities on orbit pairs, transferred by a Lusin-plus-weak-limit ``compact-set transfer lemma'' to ordered couplings built from mixed-radix digits; tails bounded without a union bound. This replaces the Brown-measure and DT-operator toolkit. The preliminary assessment is a combination bordering on new; no fresh literature search was completed for it.
\item \textbf{Affiliated transport} (Lemma 5.1): an unbounded coboundary intertwines two weights on the domain of an affiliated positive operator; invariant projections transport with equal trace via cutoffs. Incremental (Mori 2020).
\end{itemize}

\paragraph{Reused.}

Domar's convex-envelope theorems (reproved); Douglas--Pearcy equivalence of the hyperinvariant and transitive-commutant problems; Dykema 2005 on projections of hyperinvariant subspaces; Anantharaman--Popa crossed products; Borel--Cantelli; Fejér kernels.

\paragraph{Caveats.}

The repository provides two formalization routes; neither was compiled or kernel-checked for this survey. The companion's estimates ($m^{(7/10)m}$, $m^{1000}$) were not independently verified. Remark 1.5 asserts a published 2017 corollary is false. The companion estimates and comparisons with earlier equivalence claims require further checking.

\subsubsection{\texorpdfstring{297 \textperiodcentered{} Naimark's problem in ZFC}{297 · Naimark's problem in ZFC}}\label{fam:297}
\badges{Lean, 16 modules, not in catalogue \badgesep{} combination \badgesep{} earlier related result: Tanaka}

The source constructs a unital, infinite-dimensional, simple $C^*$-algebra with a faithful tracial state whose nonzero irreducible representations are all unitarily equivalent, in ZFC alone. Akemann--Weaver 2004 needed Jensen's diamond. The construction has length $\mathfrak{c}^{+}$; density is left in $[\mathfrak{c},\mathfrak{c}^+]$. Tanaka's preprint of 22 September 2026~\cite{src:tanaka} gives the same existence theorem plus density exactly $\mathfrak{c}$; the paper acknowledges this earlier result.

\paragraph{Extracted instruments.}

\begin{itemize}
\item \textbf{Selective bridge via reduced HNN adjunction along peaking sequences} (Thm 4.1). A decreasing projection sequence ``peaks to'' a state h if h is the only state with value 1 on all of them. Let $A$ be a unital $C^*$-algebra with a faithful tracial state and no finite-dimensional irreducible representations. For two inequivalent pure states, take decreasing peaking projection sequences starting at $1$, with equal traces $1=t_0>t_1>\cdots>0$ tending to zero. The source adjoins an intertwiner carrying one sequence onto the other via Ueda's reduced HNN extension (a corner of an amalgamated free product). The claim is that every old pure state extends uniquely, lifts are pure, and the only new equivalence is between the two chosen classes; simplicity and a faithful trace-preserving expectation are preserved. The proof uses norm decay on reduced words via excision, cross-excision, trace-null supports and ``moving the peak through the stable letter''. This replaces Akemann--Weaver's crossed product by a Kishimoto--Ozawa--Sakai automorphism of the whole algebra, which forced separability at each stage and hence $\aleph_{1}$ and diamond. Calderón--Farah's Question 8.1 asked for exactly such a partial crossed product. The initial extraction classifies this as a combination; a detailed comparison with Tanaka's construction remains necessary.
\item \textbf{Dependency-closed subiterations with reduced-HNN reordering} (Prop 3.4, Thm 5.5): each step records a countable set of earlier steps it depends on; adjunctions with data in a common base commute, so omitted steps can be postponed; every pure state of the big algebra is the unique extension of a pure state on a separable dependency-closed subalgebra. This supplies the ingredient Calderón--Farah \S{}8 identified. The initial extraction classifies this organization as a new combination, with a $C^*$-analogue of elementary-submodel closing-off.
\item \textbf{Diamond replacement}: length $\mathfrak{c}^{+}$, a fixed injective schedule, $\le \mathfrak{c}$ pure states per stage; no guessing principle at all. Incremental. The same argument appears to run at length $\mathfrak{c}$, which would give density $\mathfrak{c}$ as in Tanaka; this possible strengthening is not proved here or stated in the paper.
\end{itemize}

\paragraph{Reused.}

Ueda's HNN extensions and Blanchard--Dykema expectations; Akemann--Anderson--Pedersen excision; Akemann--Weaver cross-excision; Kadison transitivity; Kishimoto--Ozawa--Sakai homogeneity (only on separable subalgebras; reproved in Lean).

\paragraph{Caveats.}

The repository's formalization records differ: the catalogue omits the result, while the documentation reports full coverage. No Lean build was run here.

\subsubsection{\texorpdfstring{309 \textperiodcentered{} Kervaire invariant at the prime three}{309 · Kervaire invariant at the prime three}}\label{fam:309}
\badges{not formalized \badgesep{} combination \badgesep{} claims a new element exists}

With $b_j\in\mathrm{Ext}^{2,4\cdot3^{j+1}}$, the survivors are exactly $j\in\{0,2,3\}$: $b_1$ dies (Toda), $b_j$ dies for all $j \ge 4$, and $b_3$ survives, giving an order-three element in $\pi_{322}S_{(3)}$. This is a finite survivor list with a missing index $j=1$; the odd-primary Adams-family indexing differs from the classical prime-two Kervaire invariant convention. The result would also refute Belmont--Ray's Conjecture 0.1(3). HHR's 2011 odd-primary note explicitly left $\theta_{4}$ at $p=3$ open.

\paragraph{Extracted instruments.}

\begin{itemize}
\item \textbf{Finite Euler bound}: take the coefficient of $x^{27}$ in the strict coordinate change between two orientations of the norm $N_e^KMU$, norm and counit it to a class whose geometric fixed points vanish, so isotropy separation gives $a_\lambda^dz=0$; it maps to a unit in Borel $E_6$, so $a_\lambda^{235}=0$ and the Euler-localized homotopy fixed point spectral sequence must kill 1 by an explicit differential. This replaces HHR's Gap and Periodicity theorems and avoids any Tate-convergence claim. Combination: norms, geometric fixed points and isotropy separation are HHR and Meier--Shi--Zeng; using a coordinate-change coefficient as nilpotence witness had no matching precedent in the initial comparison.
\item \textbf{Exact cyclic coefficient model and fixed geometry}: proves $E_*\cong(W(k)[x_0..x_8,\Delta^{-1}]/(x_i+x_{i+3}+x_{i+6}))^\wedge$ with the $C_9$ action $x_i\mapsto x_{i+1}$ (the Belmont--Ray coefficient conjecture), then shows the $\sigma$-eigenspace is a strict formal module deformation of relative height 3 whose Dieudonné crystal determinant carries a sign character, excluding equivariant maps that could be differentials. This proves a stated conjecture by a combination of methods (Hill 2007 sketched formal-module deformations conditionally).
\item \textbf{Lattice forcing from the known survivor}: after $d_{33}$ the Tate pages are a finite array of fifteen labels modulo permanent unit translations; a mod-4 normalizer character and a parity count fix $\kappa=-1$; the detector of $b_2$ must be hit by $d_{93}$, forcing an outgoing differential on the whole lattice containing every $j \ge 4$ detector. The paper's key insight: a sphere class ``can acquire an incoming differential after passage to Tate localization'', so one must show outgoing, not vanishing. Combination or unclear; rests on exhaustive hand tables through $d_{313}$.
\end{itemize}

\paragraph{Reused.}

Adams--Novikov spectral sequence; Miller--Ravenel--Wilson 2-line; Toda differentials; Belmont--Ray CM detection via Lubin--Tate; Selick's criterion; Amelotte's $\theta_{3}$; the positive branch uses a truncated extended-cube construction from the stem-106 element, followed by Moore-spectrum and Adams--Novikov corrections (23,384 unknowns in the finite calculation). This produces stem $3\cdot108-2=322$; the ordinary ring cube would lie in stem $3\cdot106=318$ and is not the claimed new element.

\paragraph{Caveats.}

The claimed conclusion includes a new stable homotopy class. Its verification requires about 1,700 lines of differential calculations, with several statements specified only ``up to nonzero scalar'' entering a sign-sensitive character test. Although the paper has a reproducibility section, no code is supplied.

\subsubsection{\texorpdfstring{325 \textperiodcentered{} The complete Crouzeix conjecture}{325 · The complete Crouzeix conjecture}}\label{fam:325}
\badges{Lean, all sizes, matrix-valued, arbitrary Hilbert space \badgesep{} combination leaning new}

For every bounded operator $A$ on a nonzero complex Hilbert space, including infinite-dimensional spaces, and every matrix polynomial $P(z)=\sum_k B_kz^k$, the stated complete inequality is
\[
 \left\|\sum_k A^k\otimes B_k\right\|
 \le2\sup_{z\in W(A)}\|P(z)\|
 =2\max_{z\in\overline{W(A)}}\|P(z)\|.
\]
The constant 2 is sharp. The numerical range need not be closed in infinite dimension, which is why the supremum or its closure is used. The holomorphic extension concerns matrix-valued functions holomorphic on an open neighborhood of $\overline{W(A)}$; rational functions must have all their poles outside $\overline{W(A)}$. In this sense the closure is a complete 2-spectral set. Scalar proofs are presented in the 2026 preprints of Jin, Lorist--Schwenninger and Luo~\cite{src:jin,src:lorist,src:luo}. Lorist--Schwenninger explicitly explain why their proof and Jin's do not directly give the completely bounded version.

\paragraph{Extracted instruments.}

\begin{itemize}
\item \textbf{Ordered density comparison} (direct paper Thm 3.2, main Lemma 5.1). For a finite matrix $A$ and a matrix polynomial $F$ normalized to have boundary supremum norm one on an enclosing analytic convex domain, choose a unit top singular pair $F[A]x=\gamma y$ (the zero case is immediate). Compress the Crouzeix--Palencia positive field $R+R^*\succeq0$ by $[X\ Y]$ into a $2\times2$ block $\Delta$($\lambda$) of $m\times m$ boundary functions whose entries have Faber coefficients $P_k=X^*b_k(A)X$, $M_k=Y^*b_k(A)X$. For PSD $\Delta$, $\mathrm{tr}(\Delta J\Delta J)\ge0$ gives $2|M+L^*|^2_{HS}\le|P+P^*|^2_{HS}+|Q+Q^*|^2_{HS}$; integrate over the circle and apply Parseval. The cross term is at least $a^2=\sum|M_k|^2$; the diagonals are bounded by testing the functional on the two ordered products FG and GF (so no commutation is needed), with the Faber weight 2 giving $\le4a^2/\gamma^2$. Hence $2a^2\le8a^2/\gamma^2$ and $\gamma \le 2$ when $a>0$. If $a=0$, the cross functional vanishes, so its value on $F$ gives $\gamma=0$. The bounded-operator conclusion is obtained through finite compressions and approximation, without assuming a norm-attaining vector in infinite dimension. This replaces Crouzeix--Palencia's operator-norm triangle inequality ($1+\sqrt{2}$) and Lorist--Schwenninger's perturbation lemma, whose commutation hypothesis ``is not available after matrix amplification''. Combination leaning new: extremal pairs and the positive double layer are known; measuring the compressed block in boundary $L^{2}$(Hilbert--Schmidt) had no precedent in the literature examined.
\item \textbf{Exterior Faber coefficients and the positive kernel}: the Kövari--Pommerenke and Beckermann weight-2 bound transplanted from $H^\infty$ to matrix $L^{2}$. Incremental.
\item \textbf{Optimal metric and endpoint extremal} (main Thm 1.2): for a finite matrix and a bounded convex domain with regular real-analytic boundary strictly enclosing its numerical range, the Stein LMI attains the squared minimal condition number, with $\kappa\in[1,2]$; SDP duality plus a unitary-colligation realization gives an exact extremal; a common positive boundary representation follows from Sz.-Nagy. Combination of textbook interpolation theory.
\end{itemize}

\paragraph{Reused.}

Delyon--Delyon and Crouzeix--Palencia double-layer positivity; exterior conformal maps and Faber polynomials; Paulsen's similarity theorem; Stein 1952; Sz.-Nagy dilation; Katsnelson--Kheifets--Yuditskii colligations.

\paragraph{Caveats.}

The structural theorem applies separately to each enclosing domain strictly containing $\overline{W(A)}$; no limiting similarity on the closure is asserted. Luo's scalar argument also uses a top singular pair~\cite{src:luo}, so that ingredient alone is not evidence of novelty. The relation between its boundary certificates and the matrix-valued argument here requires closer comparison.

\subsubsection{\texorpdfstring{340 \textperiodcentered{} Counterexample to the nearby Lagrangian conjecture}{340 · Counterexample to the nearby Lagrangian conjecture}}\label{fam:340}
\badges{not formalized \badgesep{} combination}

For some large even N, $Q=S^9\times S^{N-1}$ admits a closed exact embedded Lagrangian $L\subset T^*Q$, diffeomorphic to Q, not Hamiltonian isotopic to the zero section. The paper does not address smooth isotopy or whether the projection is homotopic to a diffeomorphism. N is not explicit and L depends on a non-constructively chosen homotopy class. The comparison with Álvarez-Gavela--Igusa--Sullivan's construction concerns the passage from Legendrian to embedded Lagrangian families~\cite{src:ais}.

\paragraph{Extracted instruments.}

\begin{itemize}
\item \textbf{Momentum-separation embedding of a tube family} (Prop 4.2). From a family $g_b$, $b \in S^{9}$, of degree-2 homogeneous functions on $\mathbb R^N$ whose negative domains form a nonstandard smooth tube family, define the Morse family $F(b,v,w)=G_b(w)-R(b)\beta(|w|/T)\langle v,w\rangle$ with $R(b)=\Lambda e^{K\theta(b)}$, $\theta$ a height function. Homogeneity makes the critical locus a copy of $S^9\times S^{N-1}$; equal sphere momenta force $w-w'\parallel v$, and the base momenta separate the pair unless b is near a pole, where $g_b$ is quadratic and the critical equation has a unique solution. This converts a nontrivial tube family into an embedded exact Lagrangian, the step AIS could only do for Legendrians via wrinkling. Combination; the exponential scale over the parameter sphere is the paper's own device.
\item \textbf{One-critical-point rigidity} (Prop 3.1): a family over $S^{9}$ of functions homogeneous at infinity with exactly one nondegenerate critical point per fiber and trivial stable spherical class has trivial stable tube class, via Adams's injectivity of $J_{8}$, a family of h-cobordisms thickened into a fixed cylinder, Igusa stability (through Courte--Porcelli's bounds) and parametric multijet transversality. Incremental over Kragh's identification of one-critical-point functions with the homotopy fiber of the Hatcher--Waldhausen map.
\item \textbf{Torsion-invisible kernel class} (Prop 2.2): $\pi_{10}Wh^{diff}(*)_{(2)}$ has order 32 (Rognes) while the cokernel of $J_{9}$ has order $\le 4$, so a nonzero class exists in the kernel of the tube fibration's map to BG; degree 9 is chosen so $J_{8}$ and $J_{9}$ are injective. A 2-torsion class invisible to AIS's real-valued tube torsion. Combination of standard computations.
\end{itemize}

\paragraph{Reused.}

Waldhausen's tube theorem; Waldhausen--Jahren--Rognes stable parametrized h-cobordism theorem; Rognes 2002; Adams; Igusa; Sikorav--Chaperon generating-function transport; Bott periodicity.

\paragraph{Caveats.}

The argument depends on 2025--26 preprints by AIS, Abouzaid--Álvarez-Gavela--Courte--Kragh and Courte--Porcelli, including an informal comparison between models. Corner and stabilization conventions require detailed checking; the detection argument depends on them.

\subsubsection{\texorpdfstring{376 \textperiodcentered{} Universal computation in forced Navier--Stokes flows}{376 · Universal computation in forced Navier–Stokes flows}}\label{fam:376}
\badges{lemmas formalized; main paper not \badgesep{} incremental}

For each fixed positive computable viscosity $\nu$, an effective procedure maps a Turing machine and input to a smooth mean-zero force on $\mathbb{T}^{3}$, time-independent for $t \ge 1$, such that Navier--Stokes from rest has a global smooth solution and a marked material particle enters a fixed region in finite time iff the machine halts. Variants: $\mathbb{R}^{3}$ with compact support, force decaying faster than any power, an Eulerian velocity observer. The family contains nine closely related manuscripts.

\paragraph{Extracted instruments.}

\begin{itemize}
\item \textbf{Injective local compiler}: a three-track history-recording alphabet with a return-mark bit so that destination state and written symbol determine the incoming rule, making Moore-type target rectangles pairwise disjoint on the whole tape domain. This is an explicit implementation of Bennett 1973 history recording and is classified as incremental.
\item \textbf{Effective parking and routing isotopy}: compactly supported Hamiltonian pulses (translation and reciprocal-diagonal generators) and 3D localized curls that extract sources to parking columns, rescale, and insert at targets without collision, giving a smooth volume-preserving time-periodic field whose period map is the prescribed affine map near each rectangle. This combines standard box-moving by Hamiltonian isotopies with an effective description of the construction.
\item \textbf{Scale-nested memories under rapidly vanishing forcing}: unit-length physical steps but the n-th configuration written at digit scale $2^{-s_n}$ with displacement much smaller, so derivatives vanish faster than any power while computation continues; old records retained at coarser scales; coordinates alternate roles. The initial extraction identifies this elementary construction as the family's clearest candidate new instrument.
\end{itemize}

\paragraph{Reused.}

Moore's generalized shifts and gapped radix coding; Turing's printing problem; Leray's energy uniqueness; Helmholtz--Leray projection; heat-kernel bounds.

\paragraph{Caveats.}

For a chosen smooth divergence-free $U$ with $U(0)=0$, setting $f=U_t+(U\cdot\nabla)U-\nu\Delta U$ makes the Navier--Stokes identity hold with zero pressure. This residual construction does not automatically establish the required mean-zero condition, effective description, bounded derivatives, spatial support, eventual stationarity or time decay of the force; the properties prescribed in each version still require proof. Stationary Euler solutions are not generally Beltrami. The cited Cardona--Miranda--Peralta-Salas--Presas construction uses Beltrami fields, a special class of stationary Euler flows, to realize computation. The forced constructions therefore have a different scope from that construction and from unforced or low-mode-controlled flows. Computation is exact-real reachability with no finite-precision tolerance, as the papers concede. No fresh literature search was completed for this entry.

\section{\texorpdfstring{Verification and dependency scoreboard}{Verification and dependency scoreboard}}\label{sec:scoreboard}
The table preserves the initial extraction's classification, with source-status corrections where checked. ``Yes'' and file counts report repository artifacts or textual assessments from the initial extraction; they do not imply compilation or a complete dependency-closure check during either stage. The original labels mark 19 of the 38 rows as headline-formalized: 003, 017, 087, 102, 107, 157, 158, 196, 197, 248, 287, 362, 103, 112, 156, 221, 293, 297 and 325. Coverage inspected for 003 and 087 was partial (the zeta declaration and symmetric Mahler result respectively), so 19 cannot be read as 19 fully inspected family claims.

The original companion-only classification has five families: 143, 159, 263, 266 and 376. The rule counts a family once when its original row records a companion result or supporting lemmas but no headline coverage; including the lower-bound companion in 159 gives five, while excluding it gives four. These are counts of the original labels, not a current formalization census. In particular, the inspected snapshot also links a reciprocal-sum headline development for 159 and reports coverage for 007. The subsequent mathematical review did not include building or fully checking those developments, so their presence is recorded separately rather than silently changing what was checked.
{\footnotesize\setlength{\tabcolsep}{3pt}
\begin{longtable}{@{}>{\raggedright\arraybackslash}p{0.065\textwidth}>{\raggedright\arraybackslash}p{0.22\textwidth}>{\raggedright\arraybackslash}p{0.19\textwidth}>{\raggedright\arraybackslash}p{0.285\textwidth}>{\raggedright\arraybackslash}p{0.18\textwidth}@{}}
\caption{Verification and dependency scoreboard (initial extraction labels with selected source corrections).}\label{tab:scoreboard}\\
\toprule
\textbf{Family} & \textbf{Lean label} & \textbf{Statement faithful?} & \textbf{Depends on other OpenAI papers?} & \textbf{Uncited nearest ancestor} \\
\midrule
\endfirsthead
\multicolumn{5}{@{}l}{\emph{(continued)}}\\
\toprule
\textbf{Family} & \textbf{Lean label} & \textbf{Statement faithful?} & \textbf{Depends on other OpenAI papers?} & \textbf{Uncited nearest ancestor} \\
\midrule
\endhead
\bottomrule
\endfoot
\hyperref[fam:003]{003} & zeta partly inspected; Dirichlet/Hecke declarations also located & not inspected & no (cubic first-moment paper cited, unused) & --- \\ \addlinespace
\hyperref[fam:004]{004} & no & --- & \textbf{yes}: FM modularity at 2; pointwise 2-converse & --- \\ \addlinespace
\hyperref[fam:017]{017} & 867 files, not in catalogue & yes (eventual bound and sSup form) & no & Masser 1981; Fischler--Nakamaye 2014 \\ \addlinespace
\hyperref[fam:032]{032} & no & --- & no load-bearing (family 001 in a remark only); withdrawn papers never cited & --- \\ \addlinespace
\hyperref[fam:074]{074} & no & --- & 3D: alternative route only. 4D: cites 3D for a lemma that is really GWZ 2026 & Rios--Sawyer \\ \addlinespace
\hyperref[fam:087]{087} & original: symmetric; general/polar modules supplied, unbuilt & symmetric statement only in initial extraction & no & --- \\ \addlinespace
\hyperref[fam:102]{102} & 475 files; PCP and KMS reproved & yes, full UGC & no & --- \\ \addlinespace
\hyperref[fam:107]{107} & supplied; unbuilt here & yes, straight-line programs over $\mathbb{C}$ & no & Alman--Li 2026 \\ \addlinespace
\hyperref[fam:157]{157} & 276 files & yes, Mathlib chromatic number and clique minors & no & --- \\ \addlinespace
\hyperref[fam:158]{158} & $\approx 24{,}500$ lines & yes, unrestricted & no & --- \\ \addlinespace
\hyperref[fam:159]{159} & original: lower-bound companion; headline also linked, unbuilt & --- & no load-bearing & --- \\ \addlinespace
\hyperref[fam:196]{196} & 124 files incl. the classifying complex & yes & no & --- \\ \addlinespace
\hyperref[fam:197]{197} & 125 files; search terminates at $h \le 217$ & yes (non-soficity excluded) & type-F corollary only & --- \\ \addlinespace
\hyperref[fam:248]{248} & yes, F from scratch & yes & corollaries only & --- \\ \addlinespace
\hyperref[fam:287]{287} & 77 modules, in docs not catalogue & reads faithful & no (family 298 later, uncited) & --- \\ \addlinespace
\hyperref[fam:304]{304} & no & --- & no & --- \\ \addlinespace
\hyperref[fam:362]{362} & 207 files, continuation statement included & yes; uniqueness among global classical solutions & no & Cheverry--Ibrahim 2025 \\ \addlinespace
\midrule\multicolumn{5}{@{}l}{\textit{Second group}}\\ \midrule
\hyperref[fam:002]{002} & no & --- & \textbf{yes}: Selmer converse, two-primary BSD, Goldfeld ($\to$ Smith 2025 preprint) & --- \\ \addlinespace
\hyperref[fam:007]{007} & original: statements; docs report coverage, unbuilt & --- & no & --- \\ \addlinespace
\hyperref[fam:010]{010} & no & --- & main paper none; pro-modularity satellite depends on main + Hecke-dimension. External: Thorne Aug 2026 (unrefereed) & --- \\ \addlinespace
\hyperref[fam:040]{040} & no & --- & no (two unrefereed competing preprints cited, not used) & --- \\ \addlinespace
\hyperref[fam:103]{103} & 39 files; identifiers obfuscated & yes, multi-tape TM, $\mathrm{L}=\mathrm{RL}\wedge\mathrm{RL}=\mathrm{BPL}$ & no & --- \\ \addlinespace
\hyperref[fam:112]{112} & 179 files & yes, $\Sigma_{3}$ gate count, exact computation & no & --- \\ \addlinespace
\hyperref[fam:143]{143} & companion (quintic Liénard) only & --- & no & --- \\ \addlinespace
\hyperref[fam:156]{156} & 275 files, exact reals & yes, compact set, diam $\sqrt{2}$, no ten-cover & no (unrefereed GitHub note cited, not used) & Ky Fan-type refinements \\ \addlinespace
\hyperref[fam:193]{193} & no & --- & no (Lech paper cited as related) & Suslin--Voevodsky orbit lemma \\ \addlinespace
\hyperref[fam:221]{221} & 405 files, not in catalogue & yes, faithful statement & no & Mourrat--Panchenko; Barbier--Panchenko \\ \addlinespace
\hyperref[fam:263]{263} & companions only (607 + 428 files); headline none & companions faithful & companions depend on main paper's foundation estimates & --- \\ \addlinespace
\hyperref[fam:265]{265} & no & --- & PEPS companion depends on main & --- \\ \addlinespace
\hyperref[fam:266]{266} & companion ($N(6) \le 5$, Fourier vanishing outside Tao's $S_6$ class) only & --- & no & --- \\ \addlinespace
\hyperref[fam:268]{268} & no (computer-assisted, Python verifier) & --- & companion depends on main's intermediate data & --- \\ \addlinespace
\hyperref[fam:285]{285} & no & --- & methodological on the Kadison--Kaplansky companion; none on 196/197 & --- \\ \addlinespace
\hyperref[fam:293]{293} & 1,108-module closure, two formal routes & yes & mutual elementary lemmas with companion & --- \\ \addlinespace
\hyperref[fam:297]{297} & 16 modules, not in catalogue & yes, no hypotheses & no & Tanaka, arXiv 2609.26930 \\ \addlinespace
\hyperref[fam:309]{309} & no & --- & no & --- \\ \addlinespace
\hyperref[fam:325]{325} & 199 files, five challenges & yes, full cb statement & no & Luo 2026 (novelty comparison incomplete) \\ \addlinespace
\hyperref[fam:340]{340} & no & --- & no; external: AIS, AACK, Courte--Porcelli 2025--26 preprints & --- \\ \addlinespace
\hyperref[fam:376]{376} & 4 of 9 manuscripts, engineering lemmas & --- & dependencies among nine public manuscripts; main self-contained & Agrachev--Sarychev (not engaged) \\ \addlinespace
\end{longtable}}

Comparator challenge files contain target statements, often with \texttt{sorry} placeholders. Their JSON configurations point to separate solution modules; a placeholder in a challenge does not determine whether its solution is complete. For example, \texttt{GeneralMahler.json} points to a solution theorem distinct from its challenge, and the 003 configurations point to zeta, Dirichlet and Hecke declarations. These files were read, not compiled, and their transitive dependencies were not exhaustively inspected. Textual inspection and standard-axiom configurations do not establish kernel verification of the 19 labeled families. The configurations inspected in the subsequent mathematical review disable the optional nanoda recheck.

\section{\texorpdfstring{Limits of the extraction and mathematical review}{Limits of the extraction and mathematical review}}\label{sec:limits}

\begin{itemize}
\item \textbf{Reading is not verification.} The initial extraction and subsequent mathematical review included checks of selected statements, hypotheses and identities; those checks do not independently certify every underlying proof. No Lean formalization was built or kernel-checked during either stage. Absence of a check in this survey says nothing about whether other people have checked a manuscript.
\item \textbf{Novelty assessments depend on the search scope.} ``No precedent found'' means that the literature comparison described in the card found no matching instrument; a device in a 1980s Russian paper or an unpublished talk could be missed. The verdicts are provisional judgments relative to the literature searched; the subsequent mathematical review did not independently validate them, and mixed labels are not counted as disjoint classes.
\item \textbf{Choice of thirty-eight.} The release has 372 families; this set was selected for the prominence of the problems, which favors major conjectures and particularly strong claims. Families centered on improvements to existing bounds were not examined systematically.
\item \textbf{Search coverage.} Some entries in the second group used directly retrieved primary sources without a fresh literature search. Their novelty assessments have correspondingly narrower support.
\item \textbf{Priority and AI lineage.} At least four results have close human or AI-assisted neighbours predating the release (Naimark: Tanaka; Crouzeix: Jin, Lorist--Schwenninger, Luo; Borsuk: a 63-dimensional GitHub note; Serre positivity: a December 2025 AI note). This survey records these; it does not adjudicate priority.
\item \textbf{Scope of the extracted novelty.} The extraction identified devices within proofs. It did not identify a new foundational theory, a new category of object comparable to motives or perfectoid spaces, or a new conjecture. This describes what the extraction found in the selected families; it is not an exhaustive claim about their contents or the full release.
\end{itemize}

\section{Conclusion and future directions}\label{sec:conclusion}

This paper extracts candidate new mathematical instruments from 38 families in OpenAI's October 2026 release. The family entries describe their roles, antecedents and caveats; the maps, eleven cross-cutting patterns and dependency scoreboard make those descriptions easier to compare. The initial assessments identify new constructions and proof mechanisms, applications of established theorems in different settings, and adaptations whose earlier sources need closer acknowledgment. These findings remain subject to the limits of the literature comparisons and mathematical checks reported here.

The project begins with extraction and then asks what can be established and used. The signed-restriction case study illustrates that sequence: an instrument identified in a source paper becomes an input to a further theorem. This is the reason to investigate the release as a treasure mine. Three directions would extend the work:
\begin{enumerate}[label=(\arabic*)]
\item \textbf{Independent verification.} The cards identify where human checking would matter most: headline proofs not independently checked here (among them Hilbert's tenth problem over~$\mathbb{Q}$, the rational Hodge conjecture for CM abelian varieties, Kakeya, Erd\H{o}s's reciprocal-sum conjecture, Hilbert--Smith and the general Mahler conjecture), load-bearing dependencies on other unverified manuscripts (for example family~\hyperref[fam:004]{004} on family~\hyperref[fam:010]{010}), and specific flagged steps such as Proposition~5.14 of family~\hyperref[fam:002]{002} or the proposed extension of family~\hyperref[fam:017]{017} to $\mu(\log r)=2$ for $r\in\mathbb Q_{>0}\setminus\{1\}$, which requires a separate argument. Expert refereeing of these points, and of the novelty verdicts based on incomplete literature searches, would strengthen these assessments.
\item \textbf{Formalization.} The nineteen original headline labels include partial inspection for 003 and 087 and are not a count of fully verified results. Repository documentation and linked solution modules should first be reconciled theorem by theorem, including the additional coverage reported for 007 and 159, then compiled and kernel-checked independently. Companion proofs must be kept separate from their families' headline statements.
\item \textbf{Extending the extraction to all 372 families.} The extraction covers a sample selected for prominent problems. Applying the same protocol to the whole release, including families that improve existing bounds, would turn the catalogue of instruments into a comprehensive, searchable resource, indexed by mechanism as well as by problem, and would test whether the patterns found here hold across the release.
\end{enumerate}

\begin{samepage}
\section*{Acknowledgments}
\addcontentsline{toc}{section}{Acknowledgments}
Claude Fable 5.1 was used for the initial extraction of mathematical instruments, literature comparisons and preliminary novelty assessments. The author thanks OpenAI\textquoteright{}s GPT-6 for assistance with mathematical review and draft revision. The author takes responsibility for the establishment of this research direction, design of the methodology, the curation and interpretation of the extracted material, the development of the use case studies, and the presentation of the final manuscript.
\par\end{samepage}

\end{document}